\documentclass[a4paper,11pt]{article}
\usepackage{a4wide}
\usepackage{theorem}
\usepackage{amsmath}
\usepackage{array}
\usepackage{amssymb}
\usepackage{amsfonts}
\usepackage[french,english]{babel}
\usepackage{epsf}
\usepackage{epsfig}

\newtheorem{theo}{\indent Theorem\newline}[section]

{\theorembodyfont{\rmfamily}%
\newtheorem{rem}[theo]{\noindent Remark}}
{\theorembodyfont{\rmfamily}%
 \theoremstyle{break}%
}
\newtheorem{prop}[theo]{\indent Proposition\newline}
\newtheorem{lemma}[theo]{\indent Lemma\newline}
\newtheorem{cor}[theo]{\indent Corollary\newline}

 \def\N{{\Bbb{N}}}
\def\Z{{\Bbb{Z}}}

\def\R{{\Bbb{R}}}
\def\C{{\Bbb{C}}}

\newcommand{\ind}{\mathop{\rm ind}\nolimits}
\newcommand{\conj}{\mathop{\rm conj}\nolimits}
\newcommand{\coker}{\mathop{\rm coker}\nolimits}
\newcommand{\res}{\mathop{\rm Res}\nolimits}

\newlength{\indentation}%
\makeatletter
\newcommand\@makefntextsans[1]{%
    \parindent 0em%
    \noindent%
    \hb@xt@0em{\hss}%
    #1}
\def\footnotetextsans{%
     \@ifnextchar [\@xfootnotenextsans%
       {\@footnotetextsans}}
\def\@xfootnotenextsans[#1]{%
  \begingroup%
     \csname c@\@mpfn\endcsname #1\relax%
  \endgroup%
  \@footnotetextsans}
\long\def\@footnotetextsans#1{\insert\footins{%
    \reset@font\footnotesize%
    \interlinepenalty\interfootnotelinepenalty%
    \splittopskip\footnotesep%
    \splitmaxdepth \dp\strutbox \floatingpenalty \@MM%
    \hsize\columnwidth \@parboxrestore%
    \color@begingroup%
      \@makefntextsans{%
        \rule\z@\footnotesep\ignorespaces#1\@finalstrut\strutbox}
    \color@endgroup}}
\makeatother

\begin{document}

\cleardoublepage
\title{Invariants of real symplectic $4$-manifolds and lower bounds in real 
enumerative geometry}
\author{Jean-Yves Welschinger }
\date{}
\maketitle

\makeatletter\renewcommand{\@makefnmark}{}\makeatother
\footnotetextsans{Keywords :  Symplectic manifold, real algebraic curve, moduli space,
enumerative geometry.}
\footnotetextsans{AMS Classification : 14N10, 14P25 , 53D05 , 53D45.}

{\bf Abstract :}

We first present the construction of the moduli space of real pseudo-holomorphic curves in a
given real symplectic manifold. Then, following the approach of Gromov and Witten 
\cite{Gro,Wit,RT}, we construct invariants under deformation
of real rational symplectic $4$-manifolds. These invariants provide lower 
bounds for
the number of real rational $J$-holomorphic curves in a given homology class
passing through a given real configuration of points. 

\section*{Introduction}

Let $(X, \omega , c_X)$ be a {\it real symplectic $4$-manifold}, that is a triple made of
a $4$-manifold $X$, a symplectic form $\omega$ on $X$ and an involution $c_X$ on $X$
such that $c_X^* \omega =- \omega$, all of them being of class $C^\infty$. The fixed 
point set of $c_X$ is called {\it the real part of $X$} and is
denoted by $\R X$. It is either empty or a smooth lagrangian 
submanifold of $(X, \omega)$. Let $d \in H_2 (X ; \Z)$ be
a homology class satisfying $c_1 (X) d > 0$, where $c_1 (X)$ is the first Chern class
of the symplectic $4$-manifold $(X , \omega)$. From Corollary $1.5$ of \cite{MDSal}, we know
that the existence of such a class forces the $4$-manifold $X$ to be rational or ruled, as
soon as $d$ is not the class of an exceptional divisor. Hence, from now on, we will assume
$(X, \omega)$ to be rational. Let $x \subset X$ be a {\it real 
configuration} of points, that is a subset invariant under $c_X$, made of $c_1 (X) d - 1$
distincts points. Denote by $r$ the number of such points which are real. Let 
${\cal J}_\omega$ be the space of almost complex structures of $X$, tamed by $\omega$,
and which are of H\"older class $C^{l,\alpha}$ where $l \geq 2$ and $\alpha \in ]0,1[$ 
are fixed. This space is a contractible Banach manifold of class $C^{l,\alpha}$ (see 
\cite{Audin}, p. $42$). Denote
by $\R {\cal J}_\omega \subset {\cal J}_\omega$ the subspace consisting of those
$J \in {\cal J}_\omega$ for which $c_X$ is $J$-antiholomorphic. It is a contractible 
Banach submanifold of class $C^{l,\alpha}$ of ${\cal J}_\omega$ (see proposition 
\ref{propRJ}). 
If $J \in \R {\cal J}_\omega$ is generic enough, then there are only finitely many
$J$-holomorphic rational curves in $X$ passing through $x$ in the homology class $d$
(see Theorem \ref{theoreg}). These curves are all nodal and irreducible. The total 
number of their double points is $\delta = \frac{1}{2} (d^2 - c_1 (X) d +2)$. Let $C$ be
such a curve which is assumed to be real. We define the {\it mass} of the curve $C$ to
be the number of its real isolated double points (see \S \ref{subsectresults} for a definition). 
For every
integer $m$ ranging from $0$ to $\delta$, denote by $n_d (m)$ the total number of 
real $J$-holomorphic rational curves of mass $m$ in $X$ passing through $x$ and realizing
the homology class $d$. Then define :
$$\chi_r^d (x,J) = \sum_{m=0}^\delta (-1)^m n_d (m).$$
The main result of this paper is the following (see Theorem \ref{theoprinc}) :
\begin{theo}
\label{theointro1}
The integer $\chi_r^d (x,J)$ neither depends on the choice of $J$ nor on the choice
of $x$ (provided the number of real points in this configuration is $r$).
\end{theo}
For convenience, this integer will be denoted by $\chi_r^d$, and when $r$ does not have
the same parity as $c_1 (X) d -1$, we put $\chi_r^d$ to be $0$. We then denote by
$\chi^d (T)$ the polynomial $\sum_{r=0}^{c_1 (X) d -1} \chi_r^d T^r \in \Z [T]$. 
It follows
from Theorem \ref{theointro1} that the
function $\chi : d \in H_2 (X ; \Z) \mapsto \chi^d (T) \in \Z [T]$ only depends on
the real symplectic $4$-manifold $(X , \omega , c_X)$ and is invariant under
deformation of this real symplectic $4$-manifold. This invariant is proved to be non-trivial
for degree less or equal than five in the complex projective plane, see proposition
\ref{propchi4}. As an application of
this invariant, we obtain the following lower bounds (see Corollary \ref{corenum}) :
\begin{cor}
\label{corintro}
The integer $|\chi_r^d|$ gives a lower
bound for the total number of real rational $J$-holomorphic curves of $X$ passing through 
$x$ in the homology class $d$, independently of the choice of a generic 
$J \in \R {\cal J}_\omega$.
\end{cor}

Now, let $y = (y_1 , \dots , y_{c_1 (X) d - 2})$ be a real configuration of $c_1 (X) d - 2$
distinct points of $X$, and $s$ be the number of those which are real. We assume 
$y_{c_1 (X) d - 2}$ to be real, so that $s \neq 0$. If $J \in \R {\cal J}_\omega$ is 
generic enough, then there are 
only finitely many $J$-holomorphic rational curves in $X$ passing through $y$ in the 
homology class $d$ and having a node at $y_{c_1 (X) d - 2}$. These curves are all nodal and 
irreducible. For every integer $m$ ranging from $0$ to $\delta$, denote by $\hat{n}_d^+ (m)$ 
(resp. $\hat{n}_d^- (m)$) the total number of these curves which are real, of mass $m$ and
with a non-isolated (resp. isolated) real double point at $y_{c_1 (X) d - 2}$ (see 
\S \ref{subsectenonrel}). Define then :
$$\theta_s^d (y,J) = \sum_{m=0}^\delta (-1)^m (\hat{n}_d^+ (m) - \hat{n}_d^- (m)).$$
\begin{theo}
\label{theointro2}
The integer $\theta_s^d (y,J)$ neither depends on the choice of $J$ nor on the choice
of $y$ (provided the number of real points in this configuration is $s$).
\end{theo}
Once more, for convenience, the integer $\theta_s^d (y,J)$ will be denoted by $\theta_s^d$,
and we put $\theta_s^d = 0$ when $s$ does not have the same parity as $c_1 (X) d$. This 
invariant
makes it possible to give relations in between the coefficients of the polynomial
$\chi^d$, namely (see Theorem \ref{theorel}) :
\begin{theo}
\label{theointro3}
Let $d \in H_2 (X ; \Z)$ and $r$ be an integer between $0$ and $c_1 (X)d -3$. Then
$\chi^d_{r+2} = \chi^d_r + 2 \theta^d_{r+1}.$
\end{theo}
The text is organized as follows. The first paragraph is devoted to the construction
of the moduli space ${\cal M}_g^d (x)$ (resp. $\R {\cal M}_g^d (x)$) of pseudo-holomorphic
curves (resp. real pseudo-holomorphic curves) of genus $g$, in the homology class $d$,
and passing through the given real configuration of points $x$. The space 
$\R {\cal M}_g^d (x)$ appears to be the fixed point set of a $\Z/2\Z$-action on 
${\cal M}_g^d (x)$ induced by $c_X$. The main result of this paragraph is the theorem of
regular values (see Theorem \ref{theoreg}) which states that the set of regular values
of the Fredholm projection $\pi : {\cal M}_g^d (x) \to {\cal J}_\omega$ intersects
$\R {\cal J}_\omega$ in a dense set of the second category of $\R {\cal J}_\omega$. This 
theorem is proved for $g=0$ in this paragraph and for $g > 0$ in appendix \ref{append}.
This first paragraph is independant of the other ones and is presented in the framework of 
real symplectic manifolds of 
any dimension, since this does not require more work. The second paragraph is
devoted to the definition of the invariant $\chi$ and the proof of theorem 
\ref{theointro1}. It involves in particular many genericity arguments which are given 
in \S \ref{subsectgener}. Few computations of this invariant and applications to
real enumerative geometry are given in \S \ref{subsectresults}. For this paragraph and
the third one, we 
restrict ourselves to rational curves in real rational symplectic $4$-manifolds. Finally, the
third paragraph is devoted to the definition of the invariant $\theta$,
the statements and proofs of theorems \ref{theointro2} and \ref{theointro3} and the proof
of the non-triviality of $\chi^d$ for $d=4 , 5$ in the complex projective plane. With the exception
of this non-triviality, all these results have been announced in \cite{Wels}.\\

{\bf Acknowledgements :}

I am grateful to J.-C. Sikorav for the fruitful discussions we had on the theory of 
pseudo-holomorphic curves.

\tableofcontents

\section{Moduli space of real pseudo-holomorphic curves}
\label{sectmoduli} 

\subsection{Preliminaries}

\subsubsection{Teichm\"uller space $T_{g,m}$ and mapping class group $G$}

Let $S$ be a compact connected oriented surface of genus $g$. Let $m \in \N$ and
$z = (z_1 , \dots , z_m)$ be an ordered set of $m$ distinct points in $S$. Let 
$\tau$ be a given order two permutation of the set $\{ 1, \dots , m \}$. Denote by 
${\cal D}iff^+ (S,z)$ the group of diffeomorphisms of class $C^{k+1 , \alpha}$ of $S$, 
$k \geq 1$, which preserve the orientation of $S$ and are the identity once restricted 
to $z$. Similarly, let
${\cal D}iff (S,z)$ be the group of diffeomorphisms of class $C^{k+1 , \alpha}$ of $S$,
which fix $z$ when they preserve the orientation, or induce the permutation on $z$
associated to $\tau$ otherwise. Let ${\cal D}iff^+_0 (S,z)$ be the subgroup of 
${\cal D}iff^+ (S,z)$ consisting of diffeomorphisms isotopic to the identity, and 
${\cal J}_S$ be the space of complex structures of $S$ of class $C^{k , \alpha}$ which
are compatible with the orientation of $S$. Let $s_*$ be the morphism ${\cal D}iff (S,z) \to
\Z / 2\Z$ of kernel ${\cal D}iff^+ (S,z)$. The space ${\cal J}_S$ is a contractible Banach
manifold of class $C^{k , \alpha}$ equipped with an action of the group ${\cal D}iff (S,z)$
given by :
$$(\phi , J_S) \in {\cal D}iff (S,z) \times {\cal J}_S \mapsto s_*(\phi) (\phi^{-1})^* J_S,$$
where $(\phi^{-1})^* J_S = d\phi \circ J_S \circ d\phi^{-1}$.
Denote by $T_{g,m}$ the {\it Teichm\"uller space} ${\cal J}_S / {\cal D}iff^+_0 (S,z)$, 
it is
a finite dimensional contractible manifold. We fix a complex structure on $U_{g,m}
= S \times T_{g,m}$ so that $U_{g,m}$ is the universal curve over $T_{g,m}$, and denote
by $G^+$ (resp. $G$) the group of holomorphic (resp. holomorphic or anti-holomorphic) 
automorphisms of $U_{g,m}$. When
$(g,m) \notin \{ (0,0) , (0,1) , (0,2) , (1,0) \}$, we have $G = {\cal D}iff (S,z) /
{\cal D}iff^+_0 (S,z)$, the {\it mapping class group} of $S$. The surjective morphism
$G \to \Z / 2\Z$ of kernel $G^+$ will also be denoted by $s_*$, and we put $G^- = G
\setminus G^+$. Note that the exact sequence $1 \to G^+ \to G \to \Z / 2\Z$ splits.

\subsubsection{The manifold $\R {\cal J}_\omega$}

The real structure $c_X$ of $(X , \omega)$ induces a $\Z / 2\Z$-action on ${\cal J}_\omega$
given by $\overline{c_X}^* : J \in {\cal J}_\omega \mapsto \overline{c_X}^* (J) = 
- dc_X \circ J \circ dc_X$. Denote by $\R {\cal J}_\omega$ the fixed point set of this 
action. It consists of those $J \in \R {\cal J}_\omega$ for which $c_X$ is $J$-antiholomorphic. 
Let $J_0 \in \R {\cal J}_\omega$, the involution $\overline{c_X}^*$ induces an
involution $d_{J_0} \overline{c_X}^*$ on the tangent space $T_{J_0} {\cal J}_\omega =
L^{l,\alpha} (X , \Lambda^{0,1} X \otimes_\C TX)$, where $L^{l,\alpha} (X , \Lambda^{0,1} 
X \otimes_\C TX)$ denotes the Banach space of sections of class $C^{l,\alpha}$ of the
vector bundle $\Lambda^{0,1} X \otimes_\C TX$ over $X$. Denote by $L^{l,\alpha} (X , 
\Lambda^{0,1} X \otimes_\C TX)_{+1}$ (resp. $L^{l,\alpha} (X , 
\Lambda^{0,1} X \otimes_\C TX)_{-1}$) the eigenspace of this involution associated
to the eigenvalue $+1$ (resp. $-1$), so that $T_{J_0} {\cal J}_\omega =
L^{l,\alpha} (X , 
\Lambda^{0,1} X \otimes_\C TX)_{+1} \oplus L^{l,\alpha} (X , 
\Lambda^{0,1} X \otimes_\C TX)_{-1}$.
\begin{prop}
\label{propRJ}
The fixed point set $\R {\cal J}_\omega$ of $\overline{c_X}^*$ is a Banach submanifold 
of ${\cal J}_\omega$ of
class $C^{l,\alpha}$ which is non-empty and contractible. For every $J_0 \in 
\R {\cal J}_\omega$, the tangent space $T_{J_0} \R {\cal J}_\omega$ is
$L^{l,\alpha} (X , \Lambda^{0,1} X \otimes_\C TX)_{+1}$.
\end{prop}
Note that in particular, the decomposition $T_{J_0} {\cal J}_\omega =
L^{l,\alpha} (X , 
\Lambda^{0,1} X \otimes_\C TX)_{+1} \oplus L^{l,\alpha} (X , 
\Lambda^{0,1} X \otimes_\C TX)_{-1}$ is a direct sum of locally trivial Banach 
sub-bundles.\\

{\bf Proof:}

Let us prove first that $\R {\cal J}_\omega$ is non-empty. Let $g_X$ be a Riemannian
metric on $X$ invariant under $c_X$. Let $A \in L^\infty (X , End_\R (TX))$ be such that
for every $x \in X$, $u,v \in T_x X$, $\omega_x (u,v) = g_x (Au,v)$. Denote by $A=QJ_0$
the polar decomposition of $A$, where $J_0$ is orthogonal for $g_X$. Then $J_0 \in
{\cal J}_\omega$, $\overline{c_X}^* (A) = A$ and
since $Q = \sqrt{-A^2}$, $c_X^* (Q) = Q$. We deduce that $\overline{c_X}^* (J_0) = J_0$
and thus $J_0 \in \R {\cal J}_\omega \neq \emptyset$. Now, the Cayley-S\'evennec transform 
(see
\cite{Audin}, p. $42$) provides us with a $\Z / 2\Z$-equivariant diffeomorphism between
${\cal J}_\omega$ and
$${\cal W} = \{ W \in L^{l,\alpha} (X , End_\R (TX)) \, | \, J_0 W = -W J_0 \text{ and }
1-W^t W >>0 \},$$
where ${\cal W}$ is equipped with the involution $c_X^*$. Since the fixed point set
$\R {\cal W} = {\cal W} \cap \{ W \in L^{l,\alpha} (X , End_\R (TX)) \, | \, c_X^* (W) = 
W \}$ is a contractible Banach submanifold of class $C^{l,\alpha}$ of ${\cal W}$, the
same holds for $\R {\cal J}_\omega \subset {\cal J}_\omega$. The tangent space of
$\R {\cal J}_\omega$ is then $L^{l,\alpha} (X , \Lambda^{0,1} X \otimes_\C TX)_{+1}$.
$\square$ 

\subsection{Space of pseudo-holomorphic maps and Gromov operators}

\subsubsection{The manifolds ${\cal P}^* (x)$ and $\R {\cal P}^* (x)$}

Let $x=(x_1 , \dots , x_m)$ be an ordered set of distinct points of $X$, invariant under
$c_X$. Such a set is called a {\it real configuration} of points. Let $\tau$ be the order
two permutation of $\{ 1 , \dots , m \}$ induced by $c_X$, and $d \in H_2 (X ; \Z)$ be a
fixed homology class. Denote by 
$${\cal S}^d_g (x)= \{ u \in L^{k,p} (S,X) \, | \, u_* [S] = d \text{ and } u(z) = x \}.$$
This is a Banach manifold whose tangent space at $u \in {\cal S}^d_g (x)$ is the space
$T_u {\cal S}^d_g (x) = \{ v \in L^{k,p} (S,E_u) \, | \, v(z)=0 \},$
where $E_u$ is the bundle $u^* TX$. Let ${\cal E}$ (resp. ${\cal E}'$) be the Banach bundle
over ${\cal S}^d_g (x) \times T_{g,m} \times {\cal J}_\omega$ whose fiber over
$(u,J_S , J) \in {\cal S}^d_g (x) \times T_{g,m} \times {\cal J}_\omega$ is the space
$T_u {\cal S}^d_g (x)$ (resp. the space $L^{k-1 , p} (S , \Lambda^{0,1} S \otimes_\C E_u)$).
Let $\sigma_{\overline{\partial}}$ be the section of ${\cal E}'$ defined by
$\sigma_{\overline{\partial}} (u,J_S , J) = du + J \circ du \circ J_S$, and
${\cal P}^d_g (x)$ be the set of its zeros. This is the space of  pseudo-holomorphic maps 
from $S$ to $X$ passing through $x$. Fix some Levi-Civita connection $\nabla$ on $TX$ associated to
some Riemannian metric $g_X$ invariant under $c_X$. All the induced connections on
the bundles associated to $TX$ will also be denoted by $\nabla$, for convenience. The 
linearization of
$\sigma_{\overline{\partial}}$ at $(u,J_S , J) \in {\cal P}^d_g (x)$ is defined by
(see \cite{IShev}, formula $1.2.3$) 
$$\nabla \sigma_{\overline{\partial}} (v , \stackrel{.}{J}_S , \stackrel{.}{J}) =
Dv + J \circ du \circ \stackrel{.}{J}_S + \stackrel{.}{J} \circ du \circ J_S,$$
where $D$ is the Gromov operator defined by
$$v \in {\cal E}|_{(u,J_S , J)} \mapsto D(v) = \nabla v + J \circ \nabla v  \circ J_S +
\nabla_v J \circ du \circ J_S \in {\cal E}'|_{(u,J_S , J)}.$$
Finally, denote by ${\cal P}^* (x)$ the subspace of ${\cal P}^d_g (x)$ consisting of
non-multiple maps, that is the space of triple $(u,J_S , J)$ such that $u$ cannot be
written $u' \circ \phi$ where $\phi : S \to S'$ is a non-trivial ramified covering and
$u' : S' \to X$ is a pseudo-holomorphic map.
\begin{prop}
The space ${\cal P}^* (x)$ is a Banach manifold of class $C^{l,\alpha}$ whose tangent
space at $(u,J_S , J) \in {\cal P}^* (x)$ is the space $T_{(u,J_S , J)} {\cal P}^* (x) =
\{ (v , \stackrel{.}{J}_S , \stackrel{.}{J}) \in T_{(u,J_S , J)} ({\cal S}^d_g (x) 
\times T_{g,m} \times {\cal J}_\omega))
\, | \, \nabla_{(v , \stackrel{.}{J}_S , \stackrel{.}{J})} \sigma_{\overline{\partial}} 
 = 0 \}.$ $\square$
\end{prop}
This proposition follows from the fact that at  $(u,J_S , J) \in {\cal P}^* (x)$, the
operator $\nabla \sigma_{\overline{\partial}} : T_{(u,J_S , J)} ({\cal S}^d_g (x) 
\times T_{g,m} \times {\cal J}_\omega)) \to {\cal E}'|_{(u,J_S , J)}$ is surjective (see,
for example, \cite{Shev} Corollary $2.1.3$).

The group $G$ acts on ${\cal S}^d_g (x) \times T_{g,m} \times {\cal J}_\omega$ by 
$$\phi . (u,J_S , J) = \left\{ \begin{array}{rcl}
(u \circ \phi^{-1} , (\phi^{-1})^*J_S , J) & \text{if} & s_*(\phi)=+1, \\
(c_X \circ u \circ \phi^{-1} , \overline{(\phi^{-1})}^*J_S , \overline{c_X}^*(J)) & 
\text{if} & s_*(\phi)=-1,
\end{array} \right. $$
where $\phi \in G$ and $(u,J_S , J) \in {\cal S}^d_g (x) \times T_{g,m} \times 
{\cal J}_\omega$. Note that via the fixed identification $U_{g,m} = S
\times T_{g,m}$, $\phi$ induces some diffeomorphism of $S$ which depends on $J_S \in 
 T_{g,m}$. The action of $G$ lifts to the following actions on the bundles ${\cal E}$ 
and ${\cal E}'$ :
$$(\phi , v) \in G \times {\cal E}_{(u,J_S , J)} \mapsto \left\{ \begin{array}{rcl}
v \circ \phi^{-1} \in {\cal E}_{\phi . (u,J_S , J)} & \text{if} & s_*(\phi)=+1, \\
dc_X \circ v \circ \phi^{-1}  \in {\cal E}_{\phi . (u,J_S , J)} & 
\text{if} & s_*(\phi)=-1, 
\end{array} \right. $$
and 
$$(\phi , \alpha) \in G \times {\cal E}'_{(u,J_S , J)} \mapsto \left\{ \begin{array}{rcl}
\alpha \circ d\phi^{-1} \in {\cal E}'_{\phi . (u,J_S , J)} & \text{if} & s_*(\phi)=+1, \\
dc_X \circ \alpha \circ d\phi^{-1}  \in {\cal E}'_{\phi . (u,J_S , J)} & 
\text{if} & s_*(\phi)=-1. 
\end{array} \right. $$
The section $\sigma_{\overline{\partial}}$ is obviously $G$-equivariant for these actions.
As a consequence, the manifold ${\cal P}^* (x)$ is invariant under the action of $G$.
\begin{lemma}
\label{lemmaRP}
With the exception of the identity, only the order two elements of $G^-$ may have 
non-empty fixed 
point set in ${\cal P}^* (x)$. In particular, two such involutions have disjoint fixed
point sets. $\square$
\end{lemma}
Denote by $\R {\cal P}^* (x)$ the disjoint union of the fixed point sets of the non-trivial
elements of $G$. From Lemma \ref{lemmaRP} we know that each component of 
$\R {\cal P}^* (x)$ determines uniquely an order two element of $G^- \subset G$. This 
involution
induces bundles homomorphisms on ${\cal E}|_{\R {\cal P}^* (x)}$ and 
${\cal E}'|_{\R {\cal P}^* (x)}$. We denote by  ${\cal E}_{+1}$, ${\cal E}'_{+1}$ (resp.
${\cal E}_{-1}$, ${\cal E}'_{-1}$) the eigenspace associated to the eigenvalue $+1$
(resp. $-1$) of this homomorphism, so that ${\cal E}|_{\R {\cal P}^* (x)} = {\cal E}_{+1}
\oplus {\cal E}_{-1}$ and ${\cal E}'|_{\R {\cal P}^* (x)} = {\cal E}'_{+1}
\oplus {\cal E}'_{-1}$.
\begin{prop}
\label{propRP}
The space $\R {\cal P}^* (x)$ is a Banach submanifold of ${\cal P}^* (x)$ of class
$C^{l,\alpha}$ whose tangent space at $(u,J_S , J) \in \R {\cal P}^* (x)$ is the space
$T_{(u,J_S , J)} \R {\cal P}^* (x) = \{ (v , \stackrel{.}{J}_S , \stackrel{.}{J}) \in
{\cal E}_{+1} 
\times T_{J_S} \R T_{g,m} \times T_J \R {\cal J}_\omega
\, | \, \nabla_{(v , \stackrel{.}{J}_S , \stackrel{.}{J})} \sigma_{\overline{\partial}} 
 = 0 \}.$ 
\end{prop}
\begin{rem}
To every element $(u,J_S , J) \in \R {\cal P}^* (x)$ is associated an order two element
of $G^-$.
We denote by $\R T_{g,m}$ the fixed point set of the action of this element on $T_{g,m}$.
Hence the submanifold $\R T_{g,m}$ of $T_{g,m}$ does depend on the choice of the connected
component of $\R {\cal P}^* (x)$ and thus the notation is abusive. We will however keep
this notation for convenience.
\end{rem}

{\bf Proof:}

Let $c_S$ be an element of order two of $G^-$. From Lemma \ref{lemmaRP}, it suffices to
prove that the fixed point set of $c_S$ is a Banach submanifold of ${\cal P}^* (x)$. Note
that the action of $c_S$ on the product ${\cal S}^d_g (x) \times T_{g,m} \times 
{\cal J}_\omega$ is antiholomorphic for the almost-complex structure defined by :
$$(v , \stackrel{.}{J}_S , \stackrel{.}{J}) \in T_{(u,J_S , J)} ({\cal S}^d_g (x) 
\times T_{g,m} \times {\cal J}_\omega) \mapsto (Jv , J_S \stackrel{.}{J}_S , 
J \stackrel{.}{J})  \in T_{(u,J_S , J)} ({\cal S}^d_g (x) 
\times T_{g,m} \times {\cal J}_\omega) $$
The fixed point set of this action is a Banach manifold denoted by $\R {\cal S}^d_g (x) 
\times \R T_{g,m} \times \R {\cal J}_\omega$. The restriction of 
$\sigma_{\overline{\partial}}$ to $\R {\cal S}^d_g (x) 
\times \R T_{g,m} \times \R {\cal J}_\omega$ takes values in the sub-bundle ${\cal E}'_{+1}$
of ${\cal E}'$ associated to the eigenvalue $+1$ of $c_S$, since 
$\sigma_{\overline{\partial}}$ is $G$-equivariant. It suffices then to prove that this
restriction vanishes transversally along $\R {\cal P}^* (x) \cap (\R {\cal S}^d_g (x) 
\times \R T_{g,m} \times \R {\cal J}_\omega)$, meaning that at $(u,J_S , J) \in 
\R {\cal P}^* (x) \cap (\R {\cal S}^d_g (x) 
\times \R T_{g,m} \times \R {\cal J}_\omega)$, the operator $\nabla 
\sigma_{\overline{\partial}} : T_{(u,J_S , J)} (\R {\cal S}^d_g (x) 
\times \R T_{g,m} \times \R {\cal J}_\omega) \to {\cal E}'_{+1}|_{(u,J_S , J)}$ is
surjective. But this follows from the surjectivity of $\nabla 
\sigma_{\overline{\partial}} : T_{(u,J_S , J)} ({\cal S}^d_g (x) 
\times T_{g,m} \times {\cal J}_\omega) \to {\cal E}'|_{(u,J_S , J)}$ and the 
$G$-equivariance of $\sigma_{\overline{\partial}}$. $\square$

\subsubsection{The Gromov operators $D$ and $D_\R$}

Remember that the $\C$-linear part of the Gromov operator $D$ defined in the previous 
paragraph is some $\overline{\partial}$-operator which will be denoted by 
$\overline{\partial}$, and that its $\C$-antilinear part is some order $0$ operator 
denoted by $R$. The latter is given by the formula $R_{(u,J_S , J)} (v) = N_J 
(v, du)$
where $v \in {\cal E}$ and $N_J$ is the Nijenhuis tensor of $J$. In particular,
$R \circ du = 0$, see \cite{IShev}, Lemma $1.3.1$. Let $(u,J_S , J) \in 
{\cal P}^* (x)$, the operator $\overline{\partial}$ associated to $D$ induces a
holomorphic structure on the bundle $E_u = u^* TX$ for which the morphism $du : TS \to
E_u$ is an injective analytic bundle homomorphism (see \cite{IShev}, Lemma $1.3.1$).
Denote by $E_{u,-z} = E_u \otimes {\cal O}(-z)$, $TS_{-z} = TS \otimes {\cal O}(-z)$ and
${\cal N}_{u,-z}$ the quotient sheaf $E_{u,-z}/du(TS_{-z})$. This sheaf splits under
the form ${\cal O}(N_{u,-z}) \oplus {\cal N}_u^{sing}$, where $N_{u,-z} = N_u 
\otimes {\cal O}(-z)$, $N_u$ being the normal bundle of $u(S)$ in X, and 
${\cal N}_u^{sing} = \oplus \C_{a_i}^{n_i}$. In the latter, the sum is taken over all the 
critical points $a_i$ 
of $du$ and $\C_{a_i}^{n_i}$ denotes the skycraper sheaf of fiber $\C^{n_i}$
and support $a_i$, $n_i$ being the vanishing order of $du$ at $a_i$. The operator $D$ 
induces on the quotient an operator $D^N : L^{k,p}
(S, N_{u,-z}) \to L^{k-1,p} (S, \Lambda^{0,1} S \otimes_\C N_{u,-z})$ (see \cite{IShev},
formula $1.3.5$). Denote by $H^0_D (S , E_{u,-z})$ (resp. $H^0_D (S , N_{u,-z})$) the kernel
of the operator $D$ (resp. $D^N$), and $H^1_D (S , E_{u,-z})$ 
(resp. $H^1_D (S , N_{u,-z})$) the cokernel of this operator. We also denote by
$H^0_D (S , {\cal N}_{u,-z}) = H^0_D (S , N_{u,-z}) \oplus H^0 (S , {\cal N}_u^{sing})$.
Note that since the operators $D$ and $D^N$ are elliptic, all these spaces are finite 
dimensional and do not depend on the choice of $k,p$. They satisfy the following long
exact sequence (see \cite{Shev}, Corollary $1.5.4$) :
$$0 \to H^0 (S , TS_{-z}) \to H^0_D (S , E_{u,-z}) \to H^0_D (S , N_{u,-z}) \oplus
H^0 (S , {\cal N}_u^{sing}) \to $$
\begin{eqnarray}
\to H^1 (S , TS_{-z}) \to H^1_D (S , E_{u,-z}) \to 
H^1_D (S , N_{u,-z}) \to 0.
\end{eqnarray}
Remember finally that the dual of the operator $D^N$ is given by some operator
$D^* : L^{k,p} (S , K_S \otimes_\C  N_{u,-z}^*) \to L^{k-1,p} (S , \Lambda^{0,1} S 
\otimes_\C K_S \otimes_\C  N_{u,-z}^*)$, where $K_S = \Lambda^{1,0} S$ and $D^* = R^* -
\overline{\partial}$. Thus, the Serre duality gives isomorphisms $H^0_D (S , N_{u,-z})
\cong H^1_{D^*} (S , K_S \otimes_\C  N_{u,-z}^*)$ and $H^1_D (S , N_{u,-z}) \cong 
H^0_{D^*} (S , K_S \otimes_\C  N_{u,-z}^*)$ (see \cite{Shev}, Lemma $1.5.1$).
\begin{lemma}
The operators $D : {\cal E} \to {\cal E}'$ and $D^N$ are $G$-equivariant. $\square$
\end{lemma}
Thus, over $\R {\cal P}^* (x)$, the operator $D$ restricts to some operator
${\cal E}_{+1} \to {\cal E}'_{+1}$ which will be denoted by $D_\R$. Similarly, let 
$(u, J_S ,J)
\in \R {\cal P}^* (x)$ and $c_S$ be the associated order two element of $G^-$. Denote
by $L^{k,p} (S , N_{u,-z})_{\pm 1}$ (resp. $L^{k-1,p} (S , \Lambda^{0,1} S 
\otimes_\C  N_{u,-z})_{\pm 1}$) the eigenspace associated to the eigenvalue $\pm 1$ of
the action of $c_S$ on $L^{k,p} (S , N_{u,-z})$ (resp. 
$L^{k-1,p} (S , \Lambda^{0,1} S \otimes_\C  N_{u,-z})$). Denote then by $D_\R^N$
the operator $L^{k,p} (S , N_{u,-z})_{+ 1} \to L^{k-1,p} (S , \Lambda^{0,1} S 
\otimes_\C  N_{u,-z})_{+ 1}$ induced by $D^N$.
\begin{lemma}
\label{lemmaindDR}
The operators $D_\R$ and $D_\R^N$ are Fredholm, of indices
$\ind (D_\R) = \frac{1}{2} \ind (D) = c_1 (X) d + n(1-g)$ and
$\ind (D_\R^N) = \frac{1}{2} \ind (D^N) = c_1 (X) d + (n-3)(1-g) - (n-1)m$ where
$n= \dim_\C (X)$.
\end{lemma}

{\bf Proof:}

Fix some component of $\R {\cal P}^* (x)$ and the associated element $c_S$ of order two
of $G^-$. Remember that the decomposition of $D$ in $\C$-linear and antilinear parts
writes $\overline{\partial} +R$, where $\overline{\partial}$ and $R$ are equivariant under
the action of $c_S$ (in fact under the whole $G$). The operator $\overline{\partial}$
restricts then to some operator ${\cal E}_{+1} \to {\cal E}'_{+1}$ which remains
Fredholm, of kernel (resp. cokernel) the eigenspace $\ker_{+1} (\overline{\partial})$
(resp. $\coker_{+1} (\overline{\partial})$) associated to the eigenvalue $+1$ of the action
of $c_S$ on $\ker (\overline{\partial})$ (resp. $\coker (\overline{\partial})$).
Since $R$ is of order $0$, it follows that $D_\R$ is Fredholm of index $\ind (D_\R) = 
\ind (\overline{\partial}|_{{\cal E}_{+1}}) = \frac{1}{2} \ind (D) = c_1 (X) d + n(1-g)$.
The last equality coming from Riemann-Roch theorem and the equality before from the
fact that $\overline{\partial}$ is $\C$-linear. The same arguments applied to 
$D^N$ give the result for $D_\R^N$. $\square$\\

Denote by $H^0_D (S , E_{u,-z})_{+1}$ (resp. $H^0_D (S , N_{u,-z})_{+1}$) the kernel
of the operator $D_\R$ (resp. $D^N_\R$), and by $H^1_D (S , E_{u,-z})_{+1}$ 
(resp. $H^1_D (S , N_{u,-z})_{+1}$) the cokernel of this operator. Denote also by
$H^0_D (S , {\cal N}_{u,-z})_{+1} = H^0_D (S , N_{u,-z})_{+1} \oplus 
H^0 (S , {\cal N}_u^{sing})_{+1}$.
These spaces satisfy the following long
exact sequence :
$$0 \to H^0 (S , TS_{-z})_{+1} \to H^0_D (S , E_{u,-z})_{+1} \to H^0_D (S , N_{u,-z})_{+1}
\oplus
H^0 (S , {\cal N}_u^{sing})_{+1} \to $$
\begin{eqnarray}
\label{two}
\to H^1 (S , TS_{-z})_{+1} 
\to H^1_D (S , E_{u,-z})_{+1}
 \to H^1_D (S , N_{u,-z})_{+1} \to 0.
\end{eqnarray}
Note that $H^0_D (S , E_{u,-z})_{+1}$ (resp. $H^0_D (S , N_{u,-z})_{+1}$) coincides with
the eigenspace associated to the eigenvalue $+1$ of the action of $c_S$ on
$H^0_D (S , E_{u,-z})$ (resp. $H^0_D (S , N_{u,-z})$). Denote by 
$H^0_D (S , E_{u,-z})_{-1}$ (resp. $H^0_D (S , N_{u,-z})_{-1}$) the eigenspace associated 
to the eigenvalue $-1$.
\begin{lemma}
\label{lemmaserredual}
Serre duality provides isomorphisms :

$H^0_D (S , N_{u,-z})_{+1} \cong H^1_{D^*} (S , K_S \otimes_\C  N_{u,-z}^*)_{-1}$ and 

$H^1_D (S , N_{u,-z})_{+1} \cong 
H^0_{D^*} (S , K_S \otimes_\C  N_{u,-z}^*)_{-1}$.
\end{lemma}

{\bf Proof:}

The duality between the spaces $L^{k-1,p} (S , \Lambda^{0,1} S 
\otimes_\C  N_{u,-z})$ and $L^{k,p} (S , K_S \otimes_\C N_{u,-z}^*)$ writes
$(\psi^* , \alpha) \in L^{k,p} (S , K_S \otimes_\C N_{u,-z}^*) \times L^{k-1,p} 
(S , \Lambda^{0,1} S \otimes_\C  N_{u,-z}) \mapsto \Re e \int_S <\psi^* , \alpha>$. Now,
fix some component of $\R {\cal P}^* (x)$ and the associated element $c_S$ of order two
of $G^-$. We have 
\begin{eqnarray*}
((dc_X^*)^t \circ \psi^* \circ dc_S , dc_X \circ \alpha \circ dc_S) & = & \Re e \int_S 
<(dc_X^*)^t \circ \psi^* \circ dc_S  , dc_X \circ \alpha \circ dc_S> \\
&=& \Re e \int_S <\psi^* , \alpha> \circ dc_S \\
&=& - \Re e \int_S <\psi^* , \alpha>,
\end{eqnarray*}
since $c_S$ reverses the orientation of $S$. It follows that the spaces $H^1_D 
(S , N_{u,-z})_{\pm 1}$ and $H^0_{D^*} (S , K_S \otimes_\C  N_{u,-z}^*)_{\pm1}$
are orthogonal to each other and that the spaces $H^1_D 
(S , N_{u,-z})_{\pm 1}$ and $H^0_{D^*} (S , K_S \otimes_\C  N_{u,-z}^*)_{\mp1}$
are dual to each other. The same holds for $H^0_D 
(S , N_{u,-z})_{\pm 1}$ and $H^1_{D^*} (S , K_S \otimes_\C  N_{u,-z}^*)_{\mp1}$. $\square$

\subsection{Moduli space of pseudo-holomorphic curves}

\subsubsection{The manifolds ${\cal M}_g^d (x)$ and $\R {\cal M}_g^d (x)$ and the
projections $\pi$ and $\pi_\R$}

Denote by ${\cal M}_g^d (x)$ the quotient of ${\cal P}^* (x)$ under the action of $G^+$.
The projection $\pi : (u, J_S , J) \in {\cal P}^* (x) \mapsto J \in {\cal J}_\omega$
induces on the quotient a projection ${\cal M}_g^d (x) \to {\cal J}_\omega$ still denoted
by $\pi$. Remember the following proposition (see \cite{Shev}, Theorem $2.4.3$).
\begin{prop}
\label{proppi}
Denote by $n$ the dimension of $X$ over $\C$.

1) The space  ${\cal M}_g^d (x)$ is a Banach manifold of class $C^{l,\alpha}$ and the
projection ${\cal P}^* (x) \to {\cal M}_g^d (x)$ is a principal $G^+$-bundle.

2) The projection $\pi : {\cal M}_g^d (x) \to {\cal J}_\omega$ is Fredholm of index
$\ind_\R (\pi) = 2 (c_1 (X) d + (n-3)(1-g) - (n-1)m)$. Moreover, at $[u, J_S , J] \in
{\cal M}_g^d (x)$, the kernel of $\pi$ is isomorphic to $H^0_D 
(S , {\cal N}_{u,-z})$ and its cokernel to $H^1_D 
(S , N_{u,-z})$. $\square$
\end{prop}
The manifold ${\cal M}_g^d (x)$ is equipped with an action of the group $G / G^+ \cong
\Z / 2 \Z$. Let us denote by $\R {\cal M}_g^d (x)$ the fixed point set of this action.
The restriction of $\pi$ to $\R {\cal M}_g^d (x)$ takes value in $\R {\cal J}_\omega$.
Denote by $\pi_\R$ the induced projection $\R {\cal M}_g^d (x) \to \R {\cal J}_\omega$.
\begin{prop}
\label{proppir}
The projection $\pi_\R : \R {\cal M}_g^d (x) \to \R {\cal J}_\omega$ is Fredholm
of index $\ind_\R (\pi_\R) = c_1 (X) d + (n-3)(1-g) - (n-1)m$. Moreover, at 
$[u, J_S , J] \in
\R {\cal M}_g^d (x)$, the kernel of $\pi_\R$ is isomorphic to $H^0_D 
(S , {\cal N}_{u,-z})_{+ 1}$ and its cokernel to $H^1_D 
(S , N_{u,-z})_{+ 1}$.
\end{prop}

{\bf Proof:}

The projection $\pi$ is $\Z / 2 \Z$-equivariant. Let $[u, J_S , J] \in
\R {\cal M}_g^d (x)$ and $(u, J_S , J) \in \R {\cal P}^* (x)$ mapping to this element.
Denote by $c_S$ the associated element of order two of $G^-$. Then $Im (d_{[u, J_S , J]}
\pi_\R) = Im (d \pi ) \cap L^{l,\alpha} (X , \Lambda^{0,1} X \otimes_\C TX)_{+1}$, this
image is thus closed in $L^{l,\alpha} (X , \Lambda^{0,1} X \otimes_\C TX)_{+1} = T_J
(\R {\cal J}_\omega)$. From Proposition \ref{proppi}, we know that the cokernel of
$d_{[u, J_S , J]} \pi_\R$ is finite dimensional and isomorphic to $H^1_D 
(S , N_{u,-z})_{+ 1}$. Similarly, its kernel is finite dimensional and isomorphic to $H^0_D 
(S , {\cal N}_{u,-z})_{+ 1}$. The index formula follows from the exact sequence 
(\ref{two}), from Lemma \ref{lemmaindDR} and from the Riemann-Roch formula applied to
the bundle $TS_{-z}$, since the operator $\overline{\partial}$ on this bundle is 
$\C$-linear. $\square$

\subsubsection{The theorem of regular values}

The following theorem is the main result of this first paragraph. We only give a proof
of it in genus zero and dimension $4$ here. The general case is postponed to appendix \ref{append}.
\begin{theo}
\label{theoreg}
The set of regular values of the projection $\pi : {\cal M}_g^d (x) \to {\cal J}_\omega$
intersects $\R {\cal J}_\omega$ in a dense subset of the second category of 
$\R {\cal J}_\omega$.
\end{theo}

\begin{prop}
\label{proptransv}
The submanifold $\R {\cal J}_\omega$ of ${\cal J}_\omega$ is transversal to the
restriction of $\pi$ to ${\cal M}_g^d (x) \setminus \R {\cal M}_g^d (x)$.
\end{prop}

{\bf Proof:}

Let $J \in \R {\cal J}_\omega$ and $[u, J_S , J] \in {\cal M}_g^d (x) \setminus \R 
{\cal M}_g^d (x)$. Fix some element $(u, J_S , J) \in {\cal P}^* (x)$ lifting 
$[u, J_S , J]$. Then from Proposition \ref{proppi}, $\coker (d_{[u, J_S , J]} \pi)$ is
isomorphic to $H^1_D (S , N_{u,-z})$. Let $0 \neq \psi \in H^0_{D^*} (S , K_S
\otimes_\C N_{u,-z}^*) \cong H^1_D (S , N_{u,-z})$, it suffices to prove that there
exists $\stackrel{.}{J} \in L^{l,\alpha} (X , \Lambda^{0,1} X \otimes_\C TX)_{+1} = T_J
(\R {\cal J}_\omega)$ such that $<\psi , \stackrel{.}{J} \circ du \circ J_S> \neq 0$.
But $D^* \psi = 0$ and $D^*$ is of generalized $\overline{\partial}$-type, thus
$\psi$ vanishes only at a finite number of points (see \cite{HLS}). Since $u$ is
neither real, nor multiple, there exists an open subset $U$ of $S$, disjoint from $z
\subset S$, such that $u|_U$ is an embedding, $u(U) \cap u(S \setminus U)= \emptyset$,
$c_X (u(U)) \cap u(S) = \emptyset$ and such that $\psi$ does not vanish on $U$. Then,
there exists a section $\alpha$ of $\Lambda^{0,1} S \otimes_\C E_{u , -z}$ with support in 
$U$ such that $<\psi , \alpha> \neq 0$. Let $\stackrel{.}{J} \in L^{l,\alpha} 
(X , \Lambda^{0,1} X \otimes_\C TX)$ be a section with support in a neighborhood of 
$u(U)$ such that
$\stackrel{.}{J} \circ du \circ J_S = \alpha$. The section $\stackrel{.}{J}_\R = 
\stackrel{.}{J} + \overline{c_X}^* (\stackrel{.}{J})$ then belongs to
$L^{l,\alpha} (X , \Lambda^{0,1} X \otimes_\C TX)_{+1}$ and also satisfies
$\stackrel{.}{J}_\R \circ du \circ J_S = \alpha$, hence the result. $\square$\\

{\bf Proof of Theorem \ref{theoreg} in genus zero and dimension $4$:}

From Proposition \ref{proptransv} and the theorem of Sard-Smale (see \cite{SS}), there 
exists a dense set of the second category of $\R {\cal J}_\omega$, denoted by 
${\cal U}_1$, such that
every point of $\pi^{-1} ({\cal U}_1) \setminus \R {\cal M}_0^d (x)$ is regular for $\pi$.
Similarly, from Proposition \ref{proppir} and the theorem of Sard-Smale, the set of 
regular values of $\pi_\R$ is
a dense subset of the second category of $\R {\cal J}_\omega$ denoted by ${\cal U}_2$. Then
${\cal U} = {\cal U}_1 \cap {\cal U}_2$ is suitable. Indeed, let $J \in {\cal U}$ and
$[u, J_S , J] \in \pi_\R^{-1} (J)$. Choose some element $(u, J_S , J) \in \R {\cal P}^* (x)$
lifting it and denote by $c_S$ the associated order two element of $G^-$. By hypothesis,
$H^1_D (S , N_{u,-z})_{+1} = 0$. It suffices thus to prove that 
$H^1_D (S , N_{u,-z})_{-1} = 0$. If this would not be the case, since $S$ is rational,
we would have $H^0_D (S , N_{u,-z}) = 0$ (see \cite{HLS}, Theorem $1'$). Since $u$ is real
and $H^0 (S , {\cal N}_u^{sing})$ is carried by the cuspidal points of $u$, we see that
$\dim H^0 (S , {\cal N}_u^{sing})_{+1} = \dim H^0 (S , {\cal N}_u^{sing})_{-1} =
\frac{1}{2} \dim H^0 (S , {\cal N}_u^{sing})$. From this we would obtain
$$\ind (\pi) = 2 \dim H^0 (S , {\cal N}_u^{sing})_{+1} - \dim H^1_D (S , N_{u,-z})_{-1}
< 2 \ind (\pi_\R),$$ which contradicts Proposition \ref{proppir}. $\square$

\section{The invariant $\chi$ of real rational symplectic $4$-manifolds}

\subsection{Statements of the results}
\label{subsectresults}

Let $(X, \omega , c_X)$ be a real rational symplectic $4$-manifold and ${\cal J}_\omega$ be 
the space of almost complex structures of $X$ of class $C^{l,\alpha}$ tamed by $\omega$. 
Let $C$ be a real irreducible rational pseudo-holomorphic curve of $X$ having only ordinary nodes 
as singularities, and $d \in H_2 (X ; \Z)$ be its homology class. The total number of 
double points of $C$ is given by adjonction formula and is equal to 
$\delta = \frac{1}{2} (d^2 - c_1 (X) d +2)$. The real double points of $C$ are of two
differents natures. They are either the local intersection of two real branches, or
the local intersection of two complex conjugated branches. In the first case they are
called {\it non-isolated} and in the second case they are
called {\it isolated}.
$$\vcenter{\hbox{\input{gro1.pstex_t}}}$$
We define the {\it mass} of the curve $C$ to
be the number of its real isolated double points, it is denoted by $m(C)$. This integer
satisfies the upper and lower bounds $0 \leq m(C) \leq \delta$.
Now, let $x \subset X$ be a real configuration of $c_1 (X) d - 1$
distinct points and $r$ be the number of such points which are real. Let
$J \in \R {\cal J}_\omega$, if $J$ is generic enough, then there are only finitely many
$J$-holomorphic rational curves in $X$ passing through $x$ in the homology class $d$.
Moreover, these curves are all nodal and irreducible. For every integer $m$ ranging 
from $0$ to $\delta$, denote by $n_d (m)$ the number of these curves which are real and
of mass $m$. Then define :
$$\chi_r^d (x,J) = \sum_{m=0}^\delta (-1)^m n_d (m).$$
The main result of this paper is the following:
\begin{theo}
\label{theoprinc}
Let $(X, \omega , c_X)$ be a real rational symplectic $4$-manifold, and $d \in H_2 (X ; \Z)$.
Let $x \subset X$ be a real configuration of $c_1 (X) d - 1$
distincts points and $r$ be the cardinality of $x \cap \R X$. Finally, let 
$J \in \R {\cal J}_\omega$ be an almost complex structure generic enough, so that the
integer $\chi_r^d (x,J)$ is well defined. Then, this
integer $\chi_r^d (x,J)$ neither depends on the choice of $J$ nor on the choice
of $x$ (provided the cardinality of $x \cap \R X$ is $r$).
\end{theo}
For convenience, this integer will be denoted by $\chi_r^d$, and when $r$ does not have
the same parity as $c_1 (X) d -1$, we put $\chi_r^d$ to be $0$. We then denote by
$\chi^d (T)$ the polynomial $\sum_{r=0}^{c_1 (X) d -1} \chi_r^d T^r \in \Z [T]$. It follows
from Theorem \ref{theoprinc} that the
function $\chi : d \in H_2 (X ; \Z) \mapsto \chi^d (T) \in \Z [T]$ only depends on
the real symplectic $4$-manifold $(X , \omega , c_X)$ and is invariant under
deformation of this $4$-manifold. As an application of
this invariant, we obtain the following lower bounds:
\begin{cor}
\label{corenum}
Under the hypothesis of Theorem \ref{theoprinc}, the integer $|\chi_r^d|$ gives a lower
bound for the total number of real rational $J$-holomorphic curves of $X$ passing through 
$x$ in the homology class $d$, independently of the choice of a generic 
$J \in \R {\cal J}_\omega$. $\square$
\end{cor}
Note that this number of real curves is always bounded from above by the total number
$N_d$ of rational $J$-holomorphic curves of $X$ passing through 
$x$ in the homology class $d$, which does not depend on the choice of $J$. This number $N_d$
is a Gromov-Witten invariant of the symplectic $4$-manifold $(X , \omega)$ and was computed
by Kontsevich in \cite{Kont}. One of the main problem of real enumerative geometry is, in
this context, to know if there exists a generic real almost-complex structure $J$ so that all
these rational $J$-holomorphic curves are real. The following corollary provides a criteria
for the existence of such a structure.
\begin{cor}
Under the hypothesis of Theorem \ref{theoprinc}, assume that $\chi_r^d \geq 0$ (resp.
$\chi_r^d \leq 0$). Assume that there exists a generic 
$J \in \R {\cal J}_\omega$ such that $X$ has $\frac{1}{2} (N_d - |\chi_r^d|)$ real
$J$-holomorphic curves of odd (resp. even) mass passing through 
$x$ in the homology class $d$. Then, all of the rational $J$-holomorphic curves of $X$ 
passing through 
$x$ in the homology class $d$ are real. $\square$
\end{cor}
{\bf Examples:}

1) Let $(X, \omega , c_X)$ be the complex projective plane equipped with its standard 
symplectic
form and real structure. We denote the homology classes of the complex curves of $\C P^2$ by integers.
Then $\chi^1 (T) = 1 + T^2$, $\chi^2 (T) = T + T^3 + T^5$ and 
$\chi^3 (T) = \sum_{r=0}^4 2r T^{2r}$. The latter can be obtained computing the Euler 
caracteristic of the real part of the
blown up projective plane at the nine base points of a pencil of elliptic curves, as was
noticed by V. Kharlamov (see \cite{DgKh}, Proposition $4.7.3$ or \cite{Sot}, Theorem $3.6$).
The non-triviality of the polynomials $\chi^4 (T)$ and $\chi^5 (T)$ is proved in \S
\ref{subsectnontriv}.

2) Let $(X, \omega , c_X)$ be a real smooth cubic surface of $\C P^3$, and
$l$ be the homology class of a line. Assume that $\R X$ is homeomorphic to the blown-up
real projective plane at $2k$ points, $0 \leq k \leq 3$. Then $\chi_0^l = 2k^2 + 2k + 3$.\\

Note that when $\R X = \emptyset$, Theorem \ref{theoprinc} states that the number of real rational
$J$-holomorphic curves of $X$ passing through $x$ in the homology class $d$ does not depend on
the choices of $x$ and $J$, as it is the case for the number of complex curves.
The following question arise from Corollary \ref{corenum}. Are the lower bounds given
by Corollary \ref{corenum} sharp ? In the Example $1$, for the degree $3$ and $r=8$, the
lower bound is sharp from \cite{DgKh}, Proposition $4.7.3$. Also, is it possible to define a
similar invariant using higher genus curves in real symplectic $4$-manifolds, or in real 
symplectic 
manifolds of higher dimensions ? Note that the straightforward generalization of the integer
$\chi_r^d (x , J)$ using higher genus curves, even taking into account the coherent orientation
of the complex moduli space ${\cal M}_g^d (x , J)$, certainly does depend on $x, J$. This can 
be noticed for projective curves of degree $d \geq 4$ with one nodal point, using the same
trick as for the degree 3 curves, see Example $1$.

\subsection{Genericity arguments}
\label{subsectgener}

From now on, the real symplectic $4$-manifold $(X, \omega , c_X)$ is {\bf fixed}, so that it
will not in general be mentioned in the following statements.

Denote by $B^2$ (resp. $\overline{B}^2$) the open (resp. closed) unit disk of $\C$ and by
$j_{st}$ the standard complex structure on this disk. Similarly, denote by $B^4 (\rho)$ 
the open ball of $\C^2$ of radius $\rho$ and by
$J_{st}$ (resp. $\conj$) the standard complex (resp. real) structure on this ball.
\begin{lemma}
\label{lemmashev}
Let $(J_\lambda)_{\lambda \in ]-1 , 1[}$ be a family of almost complex structures of
class $C^{l,\alpha}$ on $B^4 (2)$ depending $C^{l-1,\alpha}$-smoothly on $\lambda$ and
satisfying $\overline{\conj}^* (J_\lambda) = J_\lambda$. Let $u_0 : B^2 \to 
B^4 (1) \subset B^4 (2)$ be a real $J_0$-holomorphic map having an isolated singularity
of order $\mu$ at $0=u_0 (0)$. Then, for every $v \in \R^2$ and every integer 
$\nu \leq 2\mu +1$, there exist $\epsilon >0$ and a family of real maps $w_\lambda
\in L^{k,p} (B^2 , \C^2)$, $\lambda \in ]-1 , 1[$, such that $w_0 = 0$, 
$\stackrel{.}{w}_0 = \frac{d}{d \lambda}|_{\lambda = 0} (w_\lambda)(0)=0$ and for every
$\lambda \in ]-\epsilon , \epsilon[$, the map $u_\lambda (t) = u_0 (t) + t^\nu 
(\lambda v + w_\lambda (t))$ is $J_\lambda$-holomorphic and real. $\square$
\end{lemma}
This is a real version of Lemma $3.1.1$ of \cite{Shev}. The proof is readily the same, it
suffices to notice that the operators $\overline{\partial}^\nu , R^\nu , D^\nu$ and
$T^\nu$ are $\Z / 2\Z$-equivariants. This proof is not reproduced here.

Denote by $\R {\cal J}_{\overline{B}^2}$ the space of complex structures of $\overline{B}^2$
compatible with the complex conjugation $\conj$.
\begin{lemma}
\label{lemmacusps}
Let $\R {\cal P}' = \{ (u, J_{\overline{B}^2} , J) \in L^{k,p} (\overline{B}^2 , X) \times \R  
{\cal J}_{\overline{B}^2} \times \R {\cal J}_\omega \, | \, du + J \circ du \circ 
J_{\overline{B}^2} = 0 \text{ and } c_X \circ u = u \circ \conj \}$, and
$\R {\cal P}'_s$ be the subspace of $\R {\cal P}'$ consisting of maps having a unique 
cuspidal point which is a real ordinary cusp interior to $B^2$. Moreover, let
$(u_\lambda , J_B^\lambda , J_\lambda )_{\lambda \in ]0,1[}$ be a path of $\R {\cal P}'$
such that $(u_0 , J_B^0 , J_0 ) \in \R {\cal P}'_s$ and $du_0$ is not injective at the
point $0 \in \overline{B}^2$. Then :

1) The space $\R {\cal P}'$ is a Banach manifold of class $C^{l,\alpha}$ and 
$\R {\cal P}'_s$ is a Banach submanifold of $\R {\cal P}'$ of codimension one.

2) The path $(u_\lambda , J_B^\lambda , J_\lambda )_{\lambda \in ]0,1[}$ is transversal to
$\R {\cal P}'$ at $\lambda = 0$ if and only if $\nabla \stackrel{.}{u}_0 (T_0 B^2)$ is
not the tangent of $u_0 (B^2)$ at the cusp $u_0 (0)$. Under this condition, there exists
$\epsilon >0$ such that for every $\lambda \in ]-\epsilon , 0[$ (resp. 
$\lambda \in ]0 , \epsilon [$), the curve $u_\lambda (\overline{B}^2)$ has a non-isolated
(resp. isolated) real double point in the neighborhood of the cusp, or vice-versa.
\end{lemma}

{\bf Proof:}

Let us start with the first part of the lemma. Remember that the Gromov operator $D$,
being elliptic and defined on the compact surface with non-empty boundary $\overline{B}^2$,
is surjective. Thus, ${\cal P}' = \{ (u, J_{\overline{B}^2}, J) \in L^{k,p} 
(\overline{B}^2 , X) 
\times {\cal J}_{\overline{B}^2} \times {\cal J}_\omega \, | \, du + J \circ du \circ 
J_{\overline{B}^2} = 0  \}$ is a Banach manifold of class $C^{l,\alpha}$ and the
projection ${\cal P}' \to {\cal J}_\omega$ is everywhere a submersion. The fact that
$\R {\cal P}'$ is a Banach submanifold  of class $C^{l,\alpha}$ of ${\cal P}'$ can be
proven in the same way as Proposition \ref{propRP}. Now, let $E$ (resp. $F$) be the vector
bundle on $\R {\cal P}' \times ]-1 , 1[$ whose fiber over $((u , J_{\overline{B}^2}, J),t)$
is the vector space $T_{u(t)} \R X$ (resp. $T^*_t ]-1 , 1[$). Denote by $\Gamma$ the
section of the bundle $F \otimes E$ defined by $\Gamma ((u , J_{\overline{B}^2}, J),t) =
d_t u$. It suffices to prove that $\Gamma$ vanishes transversally over $\R {\cal P}'_s$.
Indeed, $\R {\cal P}'_s$ is locally defined as the image of $\Gamma^{-1} (0)$ under
the projection $\R {\cal P}' \times ]-1 , 1[ \to  \R {\cal P}'$, and this projection
restricted to $\Gamma^{-1} (0)$ is an embedding (see \cite{Shev}, Lemma $3.2.5$). So let
$((u_0' , J'_{\overline{B}^2}, J'),t) \in \R {\cal P}' \times ]-1 , 1[$ be such that
$d_t u_0' = 0$, and let $v$ be an element of $T_{u_0' (t)} \R X$. There exists a real
neighborhood of $Im (u_0')$ in $X$ diffeomorphic to the ball $B^4 (2)$ of $\C^2$ via
some equivariant diffeomorphism. We can then apply Lemma \ref{lemmashev} with
$\nu = 1$ and with the constant path $J'$ of almost complex structures. Let 
$(u'_\lambda , J'_{\overline{B}^2}, J')$ be the path of $\R {\cal P}'$ given by this lemma.
We have $\Gamma ((u'_\lambda , J'_{\overline{B}^2}, J'),t) = d_t u'_\lambda$, thus
$$\nabla_{((\stackrel{.}{u}'_0 , 0 , 0) , \frac{\partial}{\partial t})} \Gamma = \nabla 
\stackrel{.}{u}'_0
+ \nabla_\frac{\partial}{\partial t} (du'_0) = v \otimes  \frac{d}{d t} +
\nabla_\frac{\partial}{\partial t} (du'_0).$$
Choosing $t$ constant in $]-1 , 1[$, we deduce the surjectivity of $\nabla \Gamma$ since
the vector $v$ has been chosen arbitrarily in $T_{u'_0 (t)} \R X$. The first part of
the lemma is proved.

Now let us prove the second part of the lemma. The kernel of the projection
$\R {\cal P}' \times ]-1 , 1[ \to  \R {\cal P}'$ is generated by vectors of the form
$((0,0,0),\frac{\partial}{\partial t})$. From what we have done, 
$\nabla_{((0 , 0 , 0) , \frac{\partial}{\partial t})} \Gamma =  
\nabla_\frac{\partial}{\partial t} (du_0)$. Since the image of 
$\nabla_\frac{\partial}{\partial t} (du_0)$ is the tangent of the curve $u_0 (]-1 , 1[)$
at the real ordinary cusp, and $\nabla_{((\stackrel{.}{u}_0 , 
\stackrel{.}{J}_{\overline{B}^2}, \stackrel{.}{J}),0)} \Gamma = \nabla \stackrel{.}{u}_0$,
we deduce the transversality condition. It remains to prove that under this 
transversality condition, the topology of the real double point of $u_\lambda$ near
the cuspidal point of $u_0$ changes when we cross the wall $\R {\cal P}'_s$ at 
$\lambda = 0$. This property is independant of the choice of the point 
$(u_1 , J^1_{\overline{B}^2}, J_1)$ of $\R {\cal P}'_s$, as soon as this point can be joined
to $(u_0 , J_{\overline{B}^2}^0, J_0)$ by a smooth path of $\R {\cal P}'_s$ transversal 
to the
projection $\R {\cal P}'_s \to \R {\cal J}_\omega$. Indeed, this follows from the fact that
the projection $\R {\cal P}'_s \to \R {\cal J}_\omega$ is a submersion. Without loss of
generality, we can assume that $(X , \omega , c_X)$ is the ball $B^4 (2)$ of $\C^2$
equipped with the standard symplectic form and complex conjugation, since the problem is
local. We will prove that $(u_0 , J_{\overline{B}^2}^0, J_0)$ can be joined to the standard
real ordinary cusp of the ball $B^4 (2) \subset \C^2$. First, we can assume that $J$ is
compatible with $\omega$. Indeed, it is easy to construct an almost complex structure
$J_0 \in \R {\cal J}_\omega$ compatible with $\omega$ and such that $u_0$ is 
$J_0$-holomorphic. Moreover, the space $\{ J \in \R {\cal J}_\omega \, | \, u_0 \text{ is }
J-\text{holomorphic} \}$ is contractible, since the Cayley-S\'evennec transform identifies 
this space
with  the space $\{ W \in End(TX) \, | \, WJ_0 = -J_0 W \, , \, 1 - W^t W >>0 \, , \,
c_X^* W = W \text{ and } T (Im u_0) \text{ is invariant under } W \}$ (see \cite{Audin}, p. $42$).
Now, from the theorem of Micallef and White \cite{MiWh}, there
exist $C^1$-diffeomorphisms $\phi$ and $\psi$ of $B^4 (2)$ and $\overline{B}^2$
respectively, such that $u_0 = \phi \circ u_{st} \circ \psi^{-1}$ where
$ u_{st} : t \in \overline{B}^2 \mapsto (t^2 , t^3) \in B^4 (2)$. Let 
$(\phi_n)_{n \in \N^*}$ (resp. $(\psi_n)_{n \in \N^*}$) be a sequence of $C^{l+1 , \alpha}$
diffeomorphisms of $B^4 (2)$ (resp. of $\overline{B}^2$) converging to
$\phi$ (resp. $\psi$) in $C^1$-topology. Let $(J_n)_{n \in \N}$ be a sequence of elements
of $\R {\cal J}_\omega$ converging to  $J_0$ and for which $\phi_n \circ u_{st} \circ 
\psi_n^{-1}$ is $J_n$-holomorphic. Such a sequence can be chosen compatible with $\omega$.
When $n$ is large enough, $(\phi_n \circ u_{st} \circ 
\psi_n^{-1} , J_B^n , J_n)$ is close enough to $(u_0 , J_B^0 , J_0)$ in $\R {\cal P}'_s$
so that the topology of the curves being on the two sides of $\R {\cal P}'_s$ in 
$\R {\cal P}'$ is the same whether they are near of $(u_0 , J_B , J_0)$ or of
$(\phi_n \circ u_{st} \circ \psi_n^{-1} , J_B^n , J_n)$. Fix such a $n$, 
there exists a smooth path $(\varphi_\tau)_{\tau \in [0,1]}$ of
$C^{l+1 , \alpha}$ diffeomorphisms of $B^4 (2)$ such that $\varphi_0 = \phi_n$ and
$\varphi_1 = Id$. We can assume that every such diffeomorphism $\varphi_\tau$ has a
constant differential preserving $\omega$ at the origin, composing them by a linear 
transformation of $\C^2$ otherwise. The path $(\varphi_\tau \circ u_{st} \circ 
\psi_n^{-1} , J_B^\tau , (\varphi_0 \circ \varphi_\tau^{-1})^* J_n)$ of $\R {\cal P}'_s$
is then transversal to the projection $\R {\cal P}' \to \R {\cal J}_\omega$ and
joins $(\varphi_0 \circ u_{st} \circ 
\psi_n^{-1}, J_B^n ,  J_n)$  to $(u_{st} \circ \psi_n^{-1} , J_B^1 , \varphi_0^* J_n)$.
Indeed, restricting ourselves to some ball of smaller radius, all these structures
$(\varphi_0 \circ \varphi_\tau^{-1})^* J_n$ are  tamed by $\omega$ since they are at
the origin. Since the space $\{ J \in \R {\cal J}_\omega \, | \, u_{st} \circ \psi_n^{-1}
\text{ is } J-\text{holomorphic} \}$ is contractible and contains the standard complex
structure $J_{st}$, it suffices to prove the result for curves of
$\R {\cal P}' \setminus \R {\cal P}'_s$ in the neighborhood of
$(u_{st} \circ \psi_n^{-1} , J'_B , J_{st})$, where 
$J'_B =(u_{st} \circ \psi_n^{-1})^*J_{st}$. Now consider the path 
$(f_\lambda )_{\lambda \in ]- \epsilon , \epsilon [}$ of $J_{st}$-holomorphic maps
$t \in \overline{B}^2 \mapsto (t^2 , t^3 + \lambda t) \in B^4 (2)$. This path is
transversal to $\R {\cal P}'_s$ at $\lambda = 0$ since $\nabla \stackrel{.}{f}_\lambda =
(0,1) \otimes \frac{d}{dt}$ and the tangent of $f_0 = u_{st}$ at the cusp is
generated by the vector $(1,0)$. Moreover, for $\lambda < 0$ (resp. $\lambda > 0$),
$f_\lambda$ has a non-isolated (resp. isolated) real double point at the parameters
$t = \pm \sqrt{-\lambda}$. The same holds for the path $f_\lambda \circ \psi_n^{-1}$, hence
the result. $\square$

\begin{prop}
\label{propcusps}
Let $(X, \omega , c_X)$ be a real rational symplectic $4$-manifold, $d \in H_2 (X ; \Z)$ and
$x$ be a real configuration of 
distinct points of $X$.

1) The subspace of $\R {\cal M}^d_0 (x)$ consisting of curves having a real ordinary cusp
(resp. a non-ordinary cusp or several cusps) is an immersed submanifold of codimension one
(resp. two).

2) The subspace of $\R {\cal M}^d_0 (x)$ consisting of curves having  a real tacnode or
a real ordinary triple point (resp. a multiple point of higher order or several such
points) is an immersed submanifold of codimension one
(resp. two) transversal to the previous one.

3) The subspace of $\R {\cal M}^d_0 (x)$ consisting of curves having a real ordinary cusp,
a real tacnode or a real ordinary triple point at some point of $x \cap \R X$ is  
an immersed submanifold of codimension two.
\end{prop}

{\bf Proof:}

To begin with, let us prove the first part of the proposition. For this purpose,
fix some component of $\R {\cal M}^d_0 (x)$ and a lift ${\cal C}$ of this component
in $\R {\cal P}^* (x)$. Denote by $c_S \in G^-$ the associated involution and by
$\R S$ the fixed point set of $c_S$ in $S$. Let $(u, J_S , J) \in {\cal C}$ be a map having
an ordinary cusp at $t \in \R S$. Fix some real neighborhood of $t$ in $S$ diffeomorphic
to $\overline{B}^2$ and a real neighborhood of $u(t)$ in $X$ diffeomorphic
to $B^4 (2)$. We deduce some ``restriction map'' $rest : {\cal C} \to \R {\cal P}'$
(see Lemma \ref{lemmacusps}). From Lemma \ref{lemmashev}, this map is transversal to 
$\R {\cal P}'_s$. Thus the subspace $rest^{-1} ( \R {\cal P}'_s) \subset {\cal C}$ made
of curves having a real cuspidal point in a neighborhood of $u(t)$ is an immersed
codimension one submanifold of ${\cal C}$. Hence, the subspace of $\R {\cal M}^d_0 (x)$ 
consisting of curves having a real ordinary cusp is an immersed submanifold of codimension 
one. Moreover, it follows from this proof that the condition to have two different cusps
is transversal, so that the subspace of these curves is an immersed submanifold of 
codimension two of $\R {\cal M}^d_0 (x)$. It remains to prove that the same holds for
curves having some cuspidal point of higher order or some non-ordinary cusp. For this,
we can assume that the cuspidal point is unique. Denote by $\R {\cal M}^d_0 (x)_s$ the 
subspace of $\R {\cal M}^d_0 (x)$ consisting of curves having a unique ordinary cusp
which is thus real. Let us fix some component of $\R {\cal M}^d_0 (x)_s$ and a lift
${\cal C}_s$ of this component in $\R {\cal P}^* (x)$. Denote as before  by $c_S \in G^-$ 
the associated involution and by
$\R S$ the fixed point set of $c_S$ in $S$. Let $(u, J_S , J) \in {\cal C}_s$ and
$t \in \R S$ be the point where $du$ vanishes. Remember that the order three jet $j^3_t (u)$
of $u$ at $t$ is well defined (see \cite{IShev}, Corollary $1.4.3$). It is a polynomial
map from $T_t \R S$
to $T_{u(t)} \R X$ whose first order term vanishes since we restrict ourselves to 
${\cal C}_s$. Denote by $\Gamma_2$ the section of the space of $2$-jets over ${\cal C}_s$
which maps $(u, J_S , J)$ to the term of order two of $j^3_t (u)$. Writing $u$
in a local chart in a neighborhood of $u(t)$, we prove that this section is smooth, of
class $C^{l,\alpha}$ (see \cite{Shev}, Lemma $3.2.3$). Moreover, it follows as before
from Lemma \ref{lemmashev} with $\nu = 2$ that this section is transversal to the zero
section. Thus, the subspace of $\R {\cal M}^d_0 (x)$ consisting of curves having
a cuspidal point of order $\geq 2$ is an immersed submanifold of codimension three. Let us
now restrict ourselves to the open set ${\cal V} \subset {\cal C}_s$ on which 
$\Gamma_2$ does not vanish. Let $(u, J_S , J) \in {\cal V}$ and $t \in \R S$ be the
cuspidal point. The term of order two of $j^3_t (u)$ defines a line in 
$T_{u(t)} \R X$, it is the tangent line of the curve $u(S)$ at $u(t)$. Denote by
$N_u(t)$ the quotient of $T_{u(t)} \R X$ by this line. Projecting the term of order three
of $j^3_t (u)$ on $N_u(t)$, we define a section $\Gamma_3$ of the bundle of $3$-jets
from $T_t \R S$ to $N_u(t)$, bundle defined over ${\cal V}$. As before, this section is
smooth of class $C^{l-1,\alpha}$ (the bundle $N_u(t)$ is only of class $C^{l-1,\alpha}$),
and it follows from Lemma \ref{lemmashev} with $\nu = 3$ that it is transversal to the 
zero section. This ends the proof of the first part of Proposition \ref{propcusps}.

Now let us prove the second part of the Proposition \ref{propcusps}. Since all the cases
are proved nearly in the same way, we only give a proof in the case of the tacnode.
Let us fix ${\cal C}$ a component
of $\R {\cal P}^* (x)$, $c_S \in G^-$ the associated involution and denote by
$\R S$ the fixed point set of $c_S$ in $S$. Let
$$\hat{\cal C} = \{ ((u, J_S , J),t_1 , t_2 ) \in {\cal C} \times \R S 
\times \R S \, | \, t_1 \neq t_2 \, , \, u(t_1) = u(t_2) \text{ and } d_{t_1} u \neq 0
\neq  d_{t_2} u \}.$$
This is  a submanifold of class $C^{l , \alpha}$ of codimension two of ${\cal C} \times
\R S \times \R S$. Denote by $N$ (resp. $F$) the vector bundle of class $C^{l-1 , \alpha}$ 
over $\hat{\cal C}$ whose fiber over $((u, J_S , J),t_1 , t_2 )$ is the quotient
space $T_{u(t_1)} \R X / d_{t_1} u (T_{t_1} \R S)$ (resp. $T^*_{t_2} \R S$). Denote by
$\Theta$ the section of the bundle $F \otimes N$ defined by $\Theta ((u, J_S , J),t_1 , 
t_2 ) = d_{t_2} u$. As before, the section $\Theta$ is smooth of class  $C^{l-1 , \alpha}$,
and from Lemma \ref{lemmashev} with $\nu = 1$, it is transversal to the zero section.
Moreover, the projection $\hat{\cal C} \to {\cal C}$ restricted to
$\Theta^{-1} (0)$ is an immersion, hence the result. The transversality of this submanifold
of codimension one of $\R {\cal M}^d_0 (x)$ with the one defined in the first part of
the proposition once more follows from Lemma \ref{lemmashev} with $\nu = 1$. The proof of
the third part of the proposition is left to the reader. $\square$\\

Denote by $\R {\cal M}^d_0 (x)_s$ the immersed submanifold of codimension one of 
$\R {\cal M}^d_0 (x)$ consisting of curves having a unique cuspidal point
which is a real ordinary cusp. Let $[u, J_S , J] \in \R {\cal M}^d_0 (x)_s$, we have
$\dim H_D^0 (S , {\cal N}_{u, -z}^{sing})_{+1} = 1$. Since $\ind (\pi_\R) = 0$, it follows
that $\dim H_D^1 (S , N_{u, -z})_{+1} \geq 1$. Thus, $S$ being rational,
 $\dim H_D^0 (S , N_{u, -z})_{+1} = 0$ (see \cite{HLS}) and $\dim H_D^1 (S , 
N_{u, -z})_{+1} = 1$. Let $\psi_u$ be a generator of the vector space
$H^0_{D^*} (S , K_S \otimes_\C N_{u, -z}^*)_{-1} = H_D^1 (S , N_{u, -z})_{+1}^*$. Remember
that since $\psi_u \neq 0$ and $D^* (\psi_u) = 0$, this section $\psi_u$ vanishes at a
finite number of points.
\begin{prop}
\label{proppsiu}
The subspace of $\R {\cal M}^d_0 (x)_s$ consisting of curves $[u, J_S , J]$ for which
the section $\psi_u$ vanishes at the unique real cuspidal point of $u$ is an immersed
submanifold of $\R {\cal M}^d_0 (x)_s$ of codimension one.
\end{prop}

{\bf Proof:}

The following proof is very analogous to the proof of Lemma $4.4.3$ of \cite{Shev}.
Fix some component of $\R {\cal M}^d_0 (x)_s$ and a lift ${\cal C}_s$ of this component
in $\R {\cal P}^* (x)$. Denote by $c_S \in G^-$ the associated involution and by
$\R S$ the fixed point set of $c_S$ in $S$. For every $(u, J_S , J) \in {\cal C}_s$, we
denote by $t_u \in \R S$ the unique point at which $du$ is not injective. Let then $F$ be
the real vector bundle of rank one on ${\cal C}_s$ whose fiber over $(u, J_S , J)$ is
the vector space $(K_S \otimes N_{u, -z}^*)_{-1}|_{t_u}$. It is a vector bundle of class
$C^{l-1 , \alpha}$. Let $\psi$ be a local section of class
$C^{l-1 , \alpha}$ of the vector bundle over ${\cal C}_s$ whose fiber over
$(u, J_S , J)$ is the vector space $H^0_{D^*} (S , K_S \otimes_\C N_{u, -z}^*)_{-1}$.
We assume that $\psi$ does not vanish and we will denote as before by $\psi_u$ the
value $\psi (u, J_S , J)$. We then deduce a section $\Gamma_\psi$ of the bundle $F$ defined
by $\Gamma_\psi (u, J_S , J) = \psi_u (t_u)$. This section $\Gamma_\psi$ is of class
$C^{l-1 , \alpha}$ and we have to prove that it vanishes transversaly. For this, we fix
a Riemannian metric on $S$ and the associated Levi-Civita connection. This connection
as well as the connection $\nabla$ on $X$ induce connections on all the associated bundles,
like $F$ for instance. For convenience, all these connections will be denoted by 
$\nabla$. Hence, suppose that $(u, J_S , J) \in {\cal C}_s$ is such that  
$\Gamma_\psi (u, J_S , J) =0$, we have to prove that $\nabla|_{(u, J_S , J)} \Gamma_\psi :
T_{(u, J_S , J)} {\cal C}_s \to F_{(u, J_S , J)}$ is surjective. Let then $j_0 \in 
(K_S \otimes N_{u, -z}^*)_{-1}|_{t_u}$, we are searching for $(v , \stackrel{.}{J}_S ,
\stackrel{.}{J}) \in T_{(u, J_S , J)} {\cal C}_s$ such that $\nabla_{(v , \stackrel{.}{J}_S ,
\stackrel{.}{J})} \Gamma_\psi = j_0$. For this, it suffices to find
$(v , \stackrel{.}{J}_S ,
\stackrel{.}{J}) \in T_{(u, J_S , J)} {\cal C}_s$ and a section $\stackrel{.}{\psi} \in
L^{k , p} (S , K_S \otimes_\C N_{u, -z}^*)_{-1}$ such that $\stackrel{.}{\psi} 
(t_u) = j_0$ and 
\begin{eqnarray}
\label{four}
\forall w \in L^{k , p} (S , N_{u, -z})_{+1}, \, <D^* \stackrel{.}{\psi} , w> +
<\psi_u , \nabla_{(v , \stackrel{.}{J}_S , \stackrel{.}{J})} D (w)> = 0.
\end{eqnarray}
It is indeed proved in \cite{Shev}, lines $(4.4.12)$ to $(4.4.15)$ that this last
relation ensures that $\nabla_{(v , \stackrel{.}{J}_S ,
\stackrel{.}{J})} \Gamma_\psi = \stackrel{.}{\psi} 
(t_u)$.

Let us start to construct the section $\stackrel{.}{\psi} $. Let $\stackrel{.}{\psi}_1 \in
L^{k , p} (S , K_S \otimes_\C N_{u, -z}^*)_{-1}$ be a local section such that
$\stackrel{.}{\psi}_1 (t_u) = j_0$ and $D^* (\stackrel{.}{\psi}_1) = 0$ in a
neighborhood of $t_u$. Such a section does exist. It suffices to solve locally the equation
$D^* ( j_0 + z\phi (z))= 0 $ where the unknown $\phi$ is defined in the neighborhood of
 $t_u$. This equation is equivalent to $(z^{-1}D^* z)(\phi (z)) = z^{-1} D^* (j_0)$. 
The operator $z^{-1}D^* z$ is equivariant under the action of $c_S$ and has, once restricted
to a
neighborhood of $t_u$, a right inverse $T$ also $\Z / 2\Z$-equivariant. Thus
$\phi = T \circ z^{-1}D^* (j_0)$ is a local solution and satisfies $d c_X \circ \phi \circ
c_S = - \phi$, which provides the existence of $\stackrel{.}{\psi}_1$. Using partition of
unity, this local section is completed to a global section $\stackrel{.}{\psi} \in
L^{k , p} (S , K_S \otimes_\C N_{u, -z}^*)_{-1}$. It remains to find
$(v , \stackrel{.}{J}_S , \stackrel{.}{J}) \in T_{(u, J_S , J)} {\cal C}_s$ such that
(\ref{four}) is satisfied. Let us search for such a vector among those for which
$v =0$, $\stackrel{.}{J}_S =0$, and $\stackrel{.}{J}=0$ along $u(S)$ and in a neighborhood 
of the cusp $u(t_u)$. Such a vetor is tangent to ${\cal C}_s$ as soon as
$\overline{c_X}^* (\stackrel{.}{J}) = \stackrel{.}{J}$. From Lemma $4.2.3$ of 
\cite{Shev}, at such a
vector $(0,0,\stackrel{.}{J})$, we have $\nabla_{(0,0,\stackrel{.}{J})} D (w) =
\nabla_w \stackrel{.}{J} \circ du \circ J_S$ and thus (\ref{four}) rewrites :
$$\forall w \in L^{k , p} (S , N_{u, -z})_{+1}, \, <D^* \stackrel{.}{\psi} , w> +
<\psi_u , \nabla_w \stackrel{.}{J} \circ du \circ J_S> = 0,$$
\begin{eqnarray}
\label{five}
\text{or} \quad - D^* \stackrel{.}{\psi} = <\psi_u , \nabla \stackrel{.}{J} \circ du \circ
J_S>.
\end{eqnarray}
Outside of a
neighborhood of $t_u$, this equation determines the value of the derivative of 
$\stackrel{.}{J}$ in the normal direction of the curve $u(S)$. After integration of
this condition, we construct a solution $\stackrel{.}{J}$ satisfying (\ref{five})
and $\overline{c_X}^* (\stackrel{.}{J})= \stackrel{.}{J}$. Hence the result. $\square$

\begin{prop}
\label{propprod}
Let $(X, \omega , c_X)$ be a real rational symplectic $4$-manifold, 
$d_1 , d_2 \in H_2 (X ; \Z)$.
and $y_1 , y_2$ be two finite disjoint subsets of $X$ invariant under $c_X$. Denote by
$$\Delta = \{ ([u_1 , J_{S_1} , J] , [u_2 , J_{S_2} , J]) \in 
\R {\cal M}_0^{d_1} (y_1) \times \R {\cal M}_0^{d_2} (y_2) \, | \, u_1 (S_1) = u_2 (S_2) \}.
$$ 
Then, the projections $\pi_\R^1 : \R {\cal M}_0^{d_1} (y_1) \to \R {\cal J}_\omega$
and $\pi_\R^2 : \R {\cal M}_0^{d_2} (y_2) \to \R {\cal J}_\omega$ are transversal outside
of $\Delta$.
\end{prop}
Note that as soon as $d_1 \neq d_2$, the diagonal $\Delta$ is empty.
\begin{cor}
\label{corprod}
Let $(X, \omega , c_X)$ be a real rational symplectic $4$-manifold, $d_1 , \dots ,
d_l \in H_2 (X ; \Z)$ and $y_1 , \dots , y_l$ be finite disjoint subsets of $X$ invariant
under $c_X$. Denote by 
$$\Delta = \{ ([u_1 , J_{S_1} , J] , \dots , [u_l , J_{S_l} , J]) \in 
\R {\cal M}_0^{d_1} (y_1) \times \dots \times \R {\cal M}_0^{d_l} (y_l) \, | \, 
\exists i \neq j \text{ for which } u_i (S_i) = u_j (S_j) \}.$$
Then the fiber product $(\Pi_{i=1}^l \R {\cal M}_0^{d_i} (y_i))\setminus \Delta$ over
$\R {\cal J}_\omega$ is a Banach manifold of class $C^{l , \alpha}$. Moreover, the 
projection $(\Pi_{i=1}^l \R {\cal M}_0^{d_i} (y_i))\setminus \Delta \to \R {\cal J}_\omega$
is Fredholm of index $c_1 (X) d - l - \# y$, where $d = \sum_{i=1}^l d_i$ and
$y = \cup_{i=1}^l y_i$. $\square$
\end{cor}

{\bf Proof of Proposition \ref{propprod}:}

Let $[u_1 , J_{S_1} , J] \in \R {\cal M}_0^{d_1} (y_1)$ and
$[u_2 , J_{S_2} , J] \in \R {\cal M}_0^{d_2} (y_2)$, so that 
$\pi_\R^1 ([u_1 , J_{S_1} , J]) = \pi_\R^2 ([u_2 , J_{S_2} , J]) = J \in 
\R {\cal J}_\omega$. Fix $(u_1 , J_{S_1} , J)$ (resp. $(u_2 , J_{S_2} , J)$) a lift of
$[u_1 , J_{S_1} , J]$ (resp. $[u_2 , J_{S_2} , J]$) in $\R {\cal P}^* (y_1)$
(resp. $\R {\cal P}^* (y_2)$), so that $u_1$ is a $J$-holomorphic map
$(S_1 , z_1) \to (X , y_1)$ (resp. $(S_2 , z_2) \to (X , y_2)$). Assume that
$u_1 (S_1) \neq u_2 (S_2)$. From Proposition \ref{proppir}, $\coker (\pi_\R^i|_{[u_i , 
J_{S_i} , J] }) \cong H_D^1 (S_i , N_{u_i , -z_i})_{+1}$. Let $0 \neq \psi_1 \in
 H_{D^*}^0 (S_1 , K_{S_1} \otimes_\C 
N_{u_1 , -z_1}^*)_{-1} \cong
H_D^1 (S_1 , N_{u_1 , -z_1})_{+1}^*$ (see Lemma \ref{lemmaserredual}). It suffices to 
prove the existence of $\stackrel{.}{J} \in L^{l , \alpha} 
(X , \Lambda^{0,1} X \otimes_\C TX)_{+1}$ such that $<\psi_1 , \stackrel{.}{J} \circ du_1
\circ J_{S_1}> \neq 0$ and $<\psi_2 , \stackrel{.}{J} \circ du_2
\circ J_{S_2}> = 0$ for every $\psi_2 \in H_{D^*}^0 (S_2 , K_{S_2} \otimes_\C 
N_{u_2 , -z_2}^*)_{-1}$. Since $D^* (\psi_1) = 0$ and $\psi_1 \neq 0$, the section
$\psi_1$ vanishes only at a finite number of points (see \cite{HLS}). Since $u_1$ is
not multiple, there exists an open set $U \subset S_1 \setminus z_1$ such that $u_1|_{U}$ is 
an embedding, $u_1(U) \cap u_1(S_1 \setminus U)= \emptyset$,
$c_X (u_1(U)) \cap u_1(U) = \emptyset$ and such that $\psi_1$ does not vanish on $U$.
Moreover, since $u_1 (S_1) \neq u_2 (S_2)$, the intersection $u_1 (S_1) \cap u_2 (S_2)$
consists only of a finite number of points and thus the open set $U$ can be chosen so
that $u_1(U) \cap u_2 (S_2) = \emptyset$.

Denote by $c_{S_1}$ the order two element of $G_1^-$ whose fixed point set in 
$\R {\cal P}^* (y_1)$ contains $(u_1 , J_{S_1} , J)$, it is in fact the
$J_{S_1}$-antiholomorphic involution of $S_1$ induced by $u_1$ and $c_X$. Let then
$\alpha_U$ be a section of $\Lambda^{0,1} S_1 \otimes E_{u_1}$ with support in $U$
such that $<\psi_1 , \alpha_U > \neq 0$. There exists $\stackrel{.}{J} \in L^{l , \alpha} 
(X , \Lambda^{0,1} X \otimes_\C TX)_{+1}$, with support in a neighborhood  of $u_1 (U)$
in $X$, such that $\stackrel{.}{J} \circ du_1 \circ J_{S_1} = \alpha_U$. Denote by
$\stackrel{.}{J}_\R = \stackrel{.}{J} -dc_X \circ \stackrel{.}{J} \circ dc_X 
\in L^{l , \alpha} 
(X , \Lambda^{0,1} X \otimes_\C TX)_{+1}$. We have
\begin{eqnarray*}
<\psi_1 , \stackrel{.}{J}_\R \circ du_1 \circ J_{S_1}> & = & <\psi_1 , 
\stackrel{.}{J} \circ du_1 \circ J_{S_1}> + <\psi_1 , \overline{c_X}^* (\stackrel{.}{J})
 \circ du_1 \circ J_{S_1}> \\
&=& <\psi_1 , \alpha_U > + <\psi_1 , dc_X \circ \stackrel{.}{J} \circ du_1 \circ J_{S_1}
\circ dc_{S_1} >\\
&=& <\psi_1 , \alpha_U > - <(dc_X)^t \circ \psi_1 \circ dc_{S_1} , \alpha_U > \\
&& (\text{changing of variables, since }c_{S_1} \text{ reverses the orientation of } S_1)\\
&=& 2<\psi_1 , \alpha_U > \neq 0.
\end{eqnarray*}
Since the support of $\stackrel{.}{J}$ is disjoint from $u_2 (S_2)$, we have
$\stackrel{.}{J}_\R \circ du_2 \circ J_{S_2} =0$ and thus $<\psi_2 , 
\stackrel{.}{J}_\R \circ du_2 \circ J_{S_2}> = 0$ for every $\psi_2 \in H_{D^*}^0 
(S_2 , K_{S_2} \otimes_\C N_{u_2 , -z_2}^*)_{-1}$. $\square$\\

Under the hypothesis of Proposition \ref{propprod}, we will denote by
$\R {\cal M}_0^{d_1 , d_2} (y_1 , y_2)$  the fiber product $(\R {\cal M}_0^{d_1} (y_1) 
\times_{\R {\cal J}_\omega} \R {\cal M}_0^{d_2} (y_2)) \setminus \Delta$.
\begin{prop}
The subspace of $\R {\cal M}_0^{d_1 , d_2} (y_1 , y_2)$ consisting of couples
$([u_1 , J_{S_1} , J] , [u_2 , J_{S_2} , J]))$ for which the union 
$u_1 (S_1) \cup u_2 (S_2)$ is not nodal or has some node at some point of
$y_1 \cup y_2$, is an immersed submanifold of codimension one.
\end{prop}

{\bf Proof:}

This subspace consists of couples for which $u_1$ or $u_2 $ is not an immersion ;
or $u_1$ and $u_2 $ are immersions, but $u_1 (S_1)$ or $u_2 (S_2)$ has some multiple
points or tacnode ; 
or $u_1 (S_1)$ and $u_2 (S_2)$ are nodal curves, but the intersection $u_1 (S_1) \cap 
u_2 (S_2)$ is not transverse ; 
or $u_1 (S_1) \cup u_2 (S_2)$ is a nodal curve, but having some node at some point of
$y_1 \cup y_2$.

In the two first cases, the result follows from Proposition \ref{propcusps}. In the two 
last cases, the proof is very much analogous to the one of cases $2$ and $3$ of Proposition
\ref{propcusps}, and mainly follows from Lemma \ref{lemmashev} with $\nu = 1$. It
is left to the reader. $\square$

\subsection{Proof of Theorem \ref{theoprinc}}
\label{proofoftheoprinc}

Let $J_0$, $J_1$ be regular values of the projections $\pi : {\cal M}^{d}_0 (x) \to
{\cal J}_\omega$ and
$\pi_i : {\cal M}^{d_i}_0 (x_i) \to
{\cal J}_\omega$, for every $d_i \in H_2 (X ; \Z)$ realized by a component of a reducible
pseudo-holomorphic curve in the class $d$, and $x_i \subset x$ a real configuration of
more than $c_1(X) d_i -1 $ points. Such values exist from Theorem \ref{theoreg}. Let
$\gamma : [0,1] \to \R {\cal J}_\omega$ be a path transversal to the projections
$\pi_\R : \R {\cal M}^{d}_0 (x) \to \R {\cal J}_\omega$ and
$\pi : ({\cal M}^{d}_0 (x) \setminus \R {\cal M}^{d}_0 (x)) \to
{\cal J}_\omega$ (see Proposition \ref{proptransv}), joining $J_0$ to $J_1$. Hence,
$\R {\cal M}_\gamma = \pi_\R^{-1} (Im (\gamma))$ is a submanifold of dimension one
of $\R {\cal M}^{d}_0 (x)$, equipped with a projection $\pi_\gamma : \R {\cal M}_\gamma
\to [0,1]$ induced by $\pi_\R$.
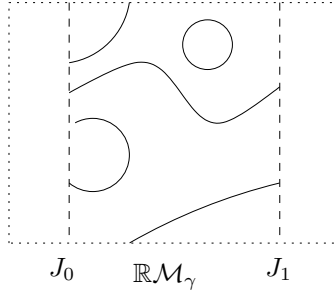
$$\vcenter{\hbox{\input{gro2.pstex_t}}}$$
The path $\gamma$ is chosen so that every element of $\R {\cal M}_\gamma$ is a nodal curve,
with the exception of a finite number of them which may have a unique real ordinary cusp,
a unique real triple point or a unique real tacnode. Moreover, this path is chosen so that
when a sequence of elements of $\R {\cal M}_\gamma$ converges in Gromov topology to a
reducible curve of $X$, then this curve has only two irreducible components, both real,
and only ordinary nodes as singularities. Finally, this path is chosen so that if
$[u , J_S , J] \in \R {\cal M}_\gamma$ has a unique real ordinary cusp at the parameter
$t_u \in \R S$, then the generator $\psi_u$ of $H^0_{D^*} (S , K_S \otimes_\C 
N_{u, -z}^*)_{-1}$ does not vanish at $t_u$. Such a choice of $\gamma$ is possible
from Propositions \ref{propcusps}, \ref{proppsiu} and \ref{propprod}.
\begin{lemma}
\label{lemmacrit}
The critical points of $\pi_\gamma$ are the curves $[u , J_S , J] \in \R {\cal M}_\gamma$
having an ordinary cusp. Moreover, all these critical points are non-degenerate.
\end{lemma}

{\bf Proof:}
From Proposition \ref{proppir}, at a point $[u, J_S , J] \in
\R {\cal M}_\gamma$, the cokernel of $d\pi_\gamma$ is isomorphic to $H^1_D 
(S , N_{u,-z})_{+ 1}$ and its kernel to
$H^0_D (S , {\cal N}_{u,-z})_{+ 1} = H^0_D (S , N_{u,-z})_{+ 1} \oplus 
H^0 (S , {\cal N}_u^{sing})_{+ 1}$. If $[u, J_S , J]$ is a critical point of $\pi_\gamma$,
then $\dim H^1_D (S , N_{u,-z})_{+ 1} = +1$. Since $S$ is rational and $D$ is of
generalized $\overline{\partial}$-type, this implies that $H^0_D (S , N_{u,-z})_{+ 1} = 0$
(see \cite{HLS}). Since $\ind (d\pi_\gamma) = 0$, we have $\dim H^0 (S , 
{\cal N}_u^{sing})_{+ 1} = 1$ and thus $u$ is not an immersion. From the hypothesis made
on $\gamma$, this implies that $u$ has a real ordinary cusp. Conversely, if
$[u, J_S , J] \in
\R {\cal M}_\gamma$ has a real ordinary cusp, then $\dim H^0 (S , 
{\cal N}_u^{sing})_{+ 1} = 1$. Thus, $\ker(d_{[u, J_S , J]} \pi_\gamma) \neq 0$ and since
$\ind (d\pi_\gamma) = 0$, $\coker(d_{[u, J_S , J]} \pi_\gamma) \neq 0$, hence the first
part of the lemma.

Now let $[u, J_S , J] \in
\R {\cal M}_\gamma$ be a critical point of $\pi_\gamma$. Fix some lift $(u, J_S , J)
\in \R {\cal P}^* (x)$ of this element, denote by $c_S \in G^-$ the associated element
of order two and by $\R S \subset S$ the fixed point set of $c_S$. Let $t_u \in \R S$ be
the point at which $du$ is not injective. We have to prove that the second order differential
$$\nabla|_{[u, J_S , J]} : H^0 (S , 
{\cal N}_u^{sing})_{+ 1} \times H^0 (S , 
{\cal N}_u^{sing})_{+ 1} \to H^1_D (S , N_{u,-z})_{+ 1}$$
is non-degenerate. Let $\psi$ be a generator of 
$H^0_{D^*} (S , K_S \otimes_\C N_{u, -z}^*)_{-1} = H_D^1 (S , N_{u, -z})_{+1}^*$ and
$(v , \stackrel{.}{J}_S)$ be a generator of $H^0 (S , 
{\cal N}_u^{sing})_{+ 1}$. Then $v$ can be written $du (\tilde{v})$ where
$\tilde{v} \in L^{k,p} (S , TS_{-z} \otimes_\C {\cal O} (t_u))_{+1}$ (see \cite{Shev}, 
Lemma $4.3.1$), that is $\tilde{v}$ is a meromorphic real vector field of $S$ having a
simple pole at $t_u$ and vanishing at $z \subset S$. From \cite{Shev}, Lemma $4.3.3$,
formula $4.3.9$, we have :
\begin{eqnarray}
<\psi , \nabla|_{[u, J_S , J]} d \pi_\gamma ((v , \stackrel{.}{J}_S) , (v , 
\stackrel{.}{J}_S))> &=& \Re e \res <\psi , \nabla d_{t_u} u (\tilde{v} , \tilde{v})> \\
&=& \Re e \lim_{\epsilon \to 0} \int_{|\xi - t_u|  = \epsilon} <\psi , \nabla d_{t_u} 
u (\tilde{v} , \tilde{v})>.
\end{eqnarray}
Now, from the computations done in the proof of Lemma $4.3.4$ of \cite{Shev} and
from Lemma $4.3.5$ of \cite{Shev}, since by hypothesis $\psi$ does not vanish at the
unique real ordinary cusp $t_u$ of $u$, the quadratic form $(6)$ is equivalent to
$w \in \R \mapsto \Re e \res_{z=0} (\frac{w^2}{z} dz)$, hence is non-degenerate.
$\square$ \\

Let $C_0$ be a real $J_0$-holomorphic nodal curve having two irreducible components
$C_1$ and $C_2$, and limit in Gromov topology of a sequence $[u_{\lambda_n} , 
J_S^{\lambda_n} , J^{\lambda_n} ]$ of elements of ${\cal M}_\gamma$, where 
$(\lambda_n)_{n \in \N}$ is a sequence of $]0,1[$ converging to some parameter
$\lambda_\infty \in ]0,1[$. Denote by $d_1 \in H_2 (X ; \Z)$ (resp. $d_2 \in H_2 (X ; \Z)$)
the homology class of $C_1$ (resp. of  $C_2$) and by $x_1 = x \cap C_1$ (resp.
$x_2 = x \cap C_2$), so that $d = d_1 + d_2$ and $x = x_1 \cup x_2$. From Propositions
\ref{propprod} and \ref{proppir} we see that, exchanging $C_1$ and $C_2$ if necessary,
we can assume that $\# (x_1) = c_1 (X) d_1 - 1$ and $\# (x_2) = c_1 (X) d_2$. Let
$[u_1 , J_{S_1} , J_0]$ and $[u_2 , J_{S_2} , J_0]$ be elements of 
$\R {\cal M}_0^{d_1} (x_1)$ and $\R {\cal M}_0^{d_2} (x_2)$ representing
$C_1$ and $C_2$ respectively. From Proposition \ref{proppir} we know that
$\dim H^1_D (S_1 , N_{u_1,-z_1})_{+ 1} \geq 0$ and $\dim H^1_D (S_2 , 
N_{u_2,-z_2})_{+ 1} \geq 1$ and from Corollary \ref{corprod} we see that these
inequalities are equalities. As a consequence, the projection $\pi_\R^1 : 
\R {\cal M}_0^{d_1} (x_1) \to \R {\cal J}_\omega$ restricts in a neighborhood of
$[u_1 , J_{S_1} , J_0]$ to a submersion on a neighborhood $V_1$ of $J_0$ in
$\R {\cal J}_\omega$. Similarly, in a neighborhood of $[u_2 , J_{S_2} , J_0]$, the projection
$\pi_\R^2 : \R {\cal M}_0^{d_2} (x_2) \to \R {\cal J}_\omega$ maps onto a codimension one
submanifold on a neighborhood $V_2$ of $J_0$ in $\R {\cal J}_\omega$. Denote by
$V = V_1 \cap V_2$ and by $H \subset V$ this codimension one submanifold. We can assume
that $V$ is connected and that  $V \setminus H$ has two connected components.
$$\vcenter{\hbox{\input{gro3.pstex_t}}}$$
Denote by $\overline{M_0^d (x)}$ (resp. $\overline{\R M_0^d (x)}$) the Gromov 
compactification of $M_0^d (x)$ (resp. $\R M_0^d (x)$), and by $\overline{\pi}$ (resp.
$\overline{\pi}_\R$) the projection $\overline{M_0^d (x)} \to {\cal J}_\omega$
(resp. $\R \overline{M_0^d (x)} \to \R {\cal J}_\omega$). Restricting $V$ if necessary, we
can assume that there exists a neighborhood $W$ of $C_0$ in $\overline{\R M_0^d (x)}$
such that $\overline{\pi} (W) = V$ and such that the image of reducible curves of $W$ under
this projection is exactly $H \subset V$. Finally, note that restricting $W$ if necessary,
we can assume that every irreducible $J$-holomorphic curve of $W$ is topologically
obtained from $C_0$ after smoothing one of the real intersection points of $C_1 \cap C_2$.
\begin{prop}
\label{propreducible}
Let $C_0$ be a real reducible $J_0$-holomorphic curve of $X$ passing through $x$ and limit
of a sequence of elements of $\R {\cal M}_\gamma$. Let $J_0 = \gamma (\lambda_0)$ for
$\lambda_0 \in ]0,1[$ and $C_1$, $C_2$ be the two irreducible components of $C_0$. 
Let $R$ be the number of real intersection points
between $C_1$ and $C_2$. Then there exist a neighborhood $W$ of $C_0$ in the Gromov 
compactification $\overline{\R {\cal M}^{d}_0 (x)}$ and $\eta > 0$ such that for every
$\lambda \in ]\lambda_0 - \eta , \lambda_0 + \eta [ \setminus \{ \lambda_0 \}$,
$\pi_\gamma^{-1} (\lambda ) \cap W$ consists exactly of $R$ real
$\gamma (\lambda)$-holomorphic curves, each of them obtained topologically by smoothing
a different real intersection point of $C_1 \cap C_2$. $\square$
\end{prop}

{\bf Proof:}

We define as above a neighborhood $W$ of $C_0$ in $\overline{\R {\cal M}^{d}_0 (x)}$, a
neighborhood $V$ of $J_0$ in $\R {\cal J}_\omega$ and a submanifold $H$ of codimension one
of $V$, such that $\overline{\pi}_\R (W) = V$ and $H$ coincide with the image under
$\overline{\pi}_\R$ of the reducible curves of $W$. Let $\eta >0$ be small enough, we will
first prove that for $\lambda \in ]\lambda_0 - \eta , \lambda_0 + \eta [ \setminus
\{ \lambda_0 \}$, $\pi_\gamma^{-1} (\lambda) \cap W$ contains at most one curve for
each real intersection point of $C_1 \cap C_2$. Otherwise, let $C'$ and $C''$ be two such
curves associated to a same intersection point of $C_1 \cap C_2$ denoted by $y_0 \in X$.
The irreducible curves $C'$ and $C''$ intersect at a finite number of points, each local
intersection being of positive multiplicity. If $W$ has been chosen small enough,
these curves have a double point in a neighborhood of each double point of $C_0$
except $y_0$. In particular, in a neighborhood of each such double point, 
these curves intersect each other in at least two points. This number of
double points is $\frac{1}{2} (d^2 - c_1 (X)d +2)$ from adjonction formula. Moreover, 
since both $C'$ and $C''$ pass through the configuration of points $x$, they intersect 
each other at each point of $x$ with multiplicity at least one. Thus, one has
$$d^2 = C' \circ C'' \geq (d^2 - c_1 (X)d +2) + c_1 (X)d - 1 = d^2 + 1,$$
which is impossible.

Let us now prove that $\pi_\gamma^{-1} (\lambda) \cap W$ actually contains exactly one
curve for each intersection point of $C_1 \cap C_2$. Let $y \in X$ be such an
intersection point. In a neighborhood of $y$, the curve $C_0$ is biholomorphic to the
standard real node ${\cal A}_0 = \{ (z^+ , z^-) \in B^2 \, | \, z^+ z^- = 0 \}$.
Following the definition $5.4.1$ of \cite{Shev}, the cylinders close to ${\cal A}_0$
are the cylinders ${\cal A}_\varphi = \{ (z^+ , z^-) \in B^2 \, | \, z^+ z^- = \varphi \}$,
for $\varphi \in B^2 (\epsilon)$, $\epsilon >0$. These cylinders form a partition of the
real analytic space ${\cal A} = \{ (z^+ , z^-) \in B^2 \, | \, |z^+ z^-| < \epsilon \}$.
Note that when the parameter $\varphi \in B^2 (\epsilon) \setminus \{ 0 \}$ is real, the cylinders
${\cal A}_\varphi$ are real and correspond topologically to the two standard ways to
smooth the real node ${\cal A}_0$, for $\varphi > 0$ and $\varphi < 0$.
$$\vcenter{\hbox{\input{gro4.pstex_t}}}$$
From Theorem $5.4.1$ of \cite{Shev} (The map $\Phi$ given in this theorem is 
$\Z/2 \Z$-equivariant for the real structures induced on ${\cal U} \times \Delta (\epsilon')$
and ${\cal P} ({\cal A})$), the real embedding of ${\cal A}_0$ in $X$ given by $C_0$
deforms into a one parameter family of real embeddings of the cylinders ${\cal A}_\varphi$,
for $\varphi \in ]-\epsilon , \epsilon[$. There exists then a continuous family 
$(J_\varphi)_{\varphi \in ]-\epsilon , \epsilon[}$ in $\R {\cal J}_\omega$ extending $J_0$,
such that for every $\varphi \in ]-\epsilon , \epsilon[ \setminus \{ 0 \}$, $J_\varphi$
differs from $J_0$ only in  a neighborhood of $\partial {\cal A}_0 \subset X$, and 
a continuous family $(C_\varphi)_{\varphi \in ]-\epsilon , \epsilon[}$ of real
$J_\varphi$-holomorphic curves extending $C_0$, such that for $\varphi \neq 0$, $C_\varphi$
is obtained topologically from $C_0$ by smoothing the real node $y$. Indeed, restricting
a little bit the real embedding of the cylinder ${\cal A}_\varphi$, the image of this
embedding can be glued to the curve $C_0 \setminus {\cal A}_0$ by adding two real small
annuli embedded in a neighborhood of $\partial {\cal A}_0$. The real curve $C_\varphi$
we thus obtain can be easily made $J_\varphi$-holomorphic for some almost-complex structure
$J_\varphi \in \R {\cal J}_\omega$ close to $J_0$ and differing from the latter only in 
a neighborhood of $\partial {\cal A}_0$.

Let us fix now $\varphi_+ \in ]0 , \epsilon[$ and $\varphi_- \in ]-\epsilon , 0[$ such that
$J_{\varphi_+}, J_{\varphi_-} \in V \setminus H$ and $C_{\varphi_+} , C_{\varphi_-} \in W$.
Deforming locally $C_{\varphi_+} , C_{\varphi_-}$ if necessary, we can assume that 
$J_{\varphi_+}, J_{\varphi_-}$ are regular values of $\pi_\R$. The almost-complex
structures $J_{\varphi_+}$ and $J_{\varphi_-}$ do not belong to the same connected
component of $V \setminus H$. Indeed, it would be otherwise possible to join them by a path
of $V \setminus H$ transversal to the projection $\pi_\R$, and there would be no obstruction
to isotop $C_{\varphi_-}$ along this path into a continuous family of $W$ to end up with
a $J_{\varphi_+}$-holomorphic curve denoted by $C'_{\varphi_-}$. The absence of such an 
obstruction follows from the fact that $C_{\varphi_-}$ can neither degenerate into a 
reducible curve nor into a cuspidal curve along this path. Since these curves are rational and 
immersed,
they are all regular points of the projection (see Proposition \ref{proppir} and \cite{HLS}).
This provides the contradiction since $C_{\varphi_+}$ and $C'_{\varphi_-}$ are two 
$J_{\varphi_+}$-holomorphic curves in $W$ obtained topologically by smoothing the same node 
of $C_0$, which is impossible from the computation done at the begining of this proof.
Now the result follows similarly. Let $\lambda_+ \in ]0 , \epsilon[$ and 
$\lambda_- \in ]-\epsilon , 0[$, the almost-complex structures $\gamma (\lambda_+)$ and
 $\gamma (\lambda_-)$ are not in the same component of $V \setminus H$. Each of them can be
joined to $J_{\varphi_+}$ or $J_{\varphi_-}$ by a path of $V \setminus H$ transversal to 
the projection $\pi_\R$. There is then no obstruction
to isotop $C_{\varphi_-}$ or $C_{\varphi_+}$ along these paths to get a 
$\gamma (\lambda_+)$ or $\gamma (\lambda_-)$-holomorphic curve which is obtained
topologically by smoothing the real node $y$ of $C_0$. Hence the result. $\square$

\begin{prop}
\label{propcusp}
Let $C_{\lambda_0} \in \R {\cal M}_\gamma$ be a critical point of $\pi_\gamma$ which is a local
maximum (resp. minimum). Then there exist a neighborhood $W$ of $C_{\lambda_0}$ in 
$\R {\cal M}_\gamma$  and $\eta > 0$ such that for every
$\lambda \in ]\lambda_0 - \eta , \lambda_0 [$ (resp. for every $\lambda \in ]\lambda_0, 
\lambda_0 + \eta [$),
$\pi_\gamma^{-1} (\lambda ) \cap W$ consists of two curves $C_\lambda^+$ and $C_\lambda^-$
satisfying $m(C_\lambda^+) = m(C_\lambda^-) + 1$, and for every $\lambda \in ]\lambda_0, 
\lambda_0 + \eta [$ (resp. for every $\lambda \in ]\lambda_0 - \eta , \lambda_0 [$),
$\pi_\gamma^{-1} (\lambda ) \cap W = \emptyset$. $\square$
\end{prop}

{\bf Proof:}

Let us assume that $C_{\lambda_0}$ is a local maximum of $\pi_\gamma$, and let us denote 
$C_{\lambda_0}$ by $[u_{\lambda_0} , J_S^{\lambda_0} , J_{\lambda_0}]$.
Since $\R {\cal M}_\gamma$ is one dimensional and $C_{\lambda_0}$ is a non-degenerate critical
point, it is clear that there exists $\eta > 0$
such that in a neighborhood $W$ of $[u_{\lambda_0} , J_S^{\lambda_0} , J_{\lambda_0}]$, 
$\pi_\gamma^{-1} (\lambda ) \cap W$ consists of two curves if 
$\lambda \in ]\lambda_0 - \eta , \lambda_0 [$, and $\pi_\gamma^{-1} (\lambda ) \cap W = 
\emptyset$ if $\lambda \in ]\lambda_0, \lambda_0 + \eta [$. The only thing to prove is that
if $\eta$ is small enough, the two curves $C_\lambda^+$ and $C_\lambda^-$ of 
$\pi_\gamma^{-1} (\lambda ) \cap W$ satisfy $m(C_\lambda^+) = m(C_\lambda^-) + 1$. From
the choice of $\gamma$ we know that the only singularities of $C_{\lambda_0}$ are nodes and
a unique real ordinary cusp. If $\eta$ is small enough, the two curves $C_\lambda^+$ and 
$C_\lambda^-$ are close enough to $C_{\lambda_0}$ so that they have a node in a neighborhood
of each node of $C_{\lambda_0}$ plus a node in a neighborhood of the cusp of $C_{\lambda_0}$.
Since the nodes which are close to nodes of $C_{\lambda_0}$ are of the same nature, we have
to prove that for one of the curves $C_\lambda^+$ or 
$C_\lambda^-$, the real node close to the cusp of $C_{\lambda_0}$ is non-isolated, and 
for the other one, it is isolated. The result indeed follows then from the definition of 
the mass.

So let us fix a parametrization $\mu \in ]-\epsilon , \epsilon[ \mapsto C_\mu \in 
\R {\cal M}_\gamma \cap W$, such that $C_0 = C_{\lambda_0}$ and $\pi_\gamma (C_\mu) = 
\pi_\gamma (C_{- \mu})$. Considering the restriction of these curves to a neighborhood
of the cusp of $C_0$ diffeomorphic to the ball $B^4 (2)$ of $\C^2$, we deduce, with the
notations of Lemma \ref{lemmacusps}, a path $(C'_\mu)_{\mu \in ]-\epsilon , \epsilon[}$
of $\R {\cal P}'$ such that $C'_0 \in \R {\cal P}'_s$. We have to prove that this path
is transversal to $\R {\cal P}'_s$ at $\mu = 0$. From the hypothesis, we know that
$v_{\lambda_0} = \frac{d}{d\mu} (C_\mu)|_{\mu = 0} \in 
H_D^0 (S , {\cal N}_{u , -z}^{sing} )_{+1}$. Denoting by $t_0 \in \R S$ the point at which
$d u_{\lambda_0}$ is not injective, we deduce from Lemma $4.3.1$ of \cite{Shev} that
$v_{\lambda_0} = d u_{\lambda_0} (w_{\lambda_0})$ for a section $w_{\lambda_0}$ of
class $L^{k-1,p}$ of the bundle $T S_{-z} \otimes_\C {\cal O} (t_0)$, that is for a vector
field of $S$ vanishing at $z$ and with a simple pole at $t_0$. Restricting ourselves to
the ball $B^4 (2)$ defined above, we can write $w_{\lambda_0} = \frac{1}{t} w'_{\lambda_0}$.
Moreover, from Corollary $3.1.3$ of \cite{IShev}, the diffeomorphism onto this ball
can be chosen in order that $u_{\lambda_0}$ writes $t \mapsto (t^2 , t^3 u'_{\lambda_0})$
with $ u'_{\lambda_0} \in L^{k,p} (B^4 (2) , \C^2)$ and $u'_{\lambda_0} (0) \neq 0$.
Thus $v_{\lambda_0} = (2 w'_{\lambda_0} , 3t w'_{\lambda_0} u'_{\lambda_0}  + t^2
w'_{\lambda_0} \frac{d}{dt} u'_{\lambda_0})$ and 
$\nabla|_{t=0} v_{\lambda_0} = (2 \nabla|_{t=0} w'_{\lambda_0} , 3 w'_{\lambda_0} 
u'_{\lambda_0} \frac{d}{dt})$. Hence $Im ( \nabla|_{t=0} v_{\lambda_0})$ is not the tangent
of $C_{\lambda_0}$ at the cusp and the transversality condition of Lemma \ref{lemmacusps}
is satisfied, which proves the result. $\square$\\

{\bf Proof of Theorem \ref{theoprinc} :}

Let $J_0 , J_1 \in \R {\cal J}_\omega$ be two regular values of the projection
$\pi : {\cal M}^d_0 (x) \to {\cal J}_\omega$, and $\gamma : [0,1] \to \R {\cal J}_\omega$
be the path fixed at the begining of \S \ref{proofoftheoprinc} joining them.
The integer $\chi_r^d (x, \gamma (\lambda))$ is then well defined for all $\lambda \in 
[0,1]$ but a finite number of values $0 < \lambda_1 < \dots < \lambda_j < 1$ corresponding
either to reducible curves, to cuspidal curves, or to curves having a real triple point
or tacnode. Since the function $\lambda \mapsto \chi_r^d (x, \gamma (\lambda))$ is
obviously constant between these values, we just have to prove that for $i \in \{ 1,
\dots , j \}$, $\chi_r^d (x, \gamma (\lambda_i^-)) = \chi_r^d (x, \gamma (\lambda_i^+))$
where $\lambda_i^-$ (resp. $\lambda_i^+$) is the left limit (resp. right limit) of 
$\lambda$ at $\lambda_i$. If $\lambda_i$ corresponds to a curve having a real triple point
or tacnode, it is straightforward and illustrated by the following pictures.

$$\vcenter{\hbox{\input{gro5.pstex_t}}}$$
\vspace{0.5cm}
$$\vcenter{\hbox{\input{gro6.pstex_t}}}$$
\vspace{0.5cm}

If $\lambda_i$ corresponds to a reducible curve, it follows from Proposition 
\ref{propreducible}, and if $\lambda_i$ corresponds to a cuspidal curve, it follows from 
Proposition \ref{propcusp}. Hence, the integer $\chi_r^d (x,J)$ does not depend on the
choice of $J \in \R {\cal J}_\omega$. Note that in contrast to the previous cases, the coefficient
$(-1)^m$ in the definition of $\chi_r^d (x,J)$ plays in this last case a crucial r\^ole to get 
the invariance. This integer $\chi_r^d (x,J)$ also does not depend on the choice of $x$,
since the group of equivariant diffeomorphisms of $X$ acts transitively on these real
configurations of points. Theorem \ref{theoprinc} is thus proved. $\square$

\section{Further study of the polynomial $\chi^d (T)$}

In the first three subparagraphs, we will give the relations between the coefficients
of the polynomial $\chi^d (T)$ in term of a new invariant $\theta$. In the last
subparagraph, we will prove the non-triviality of this polynomial in degrees $4$ and $5$
in $(\C P^2 , \omega_{std} , \conj)$ (for degree less than four, it has already been
computed in the first example given in \S \ref{subsectresults}).

\subsection{The invariant $\theta$}
\label{subsectenonrel}

Let $y = (y_1 , \dots , y_{c_1 (X) d - 2})$ be a real configuration of $c_1 (X) d - 2$
distinct points of $X$, and $s$ be the number of those which are real. 
We assume that $y_{c_1 (X) d - 2}$ is real, so that $s$ does not vanish. Let 
$J \in \R {\cal J}_\omega$ be generic enough. Then there are 
only finitely many $J$-holomorphic rational curves in $X$ in the 
homology class $d$ passing through $y$ and having an ordinary node at $y_{c_1 (X) d - 2}$. 
These curves are all nodal and 
irreducible. For every integer $m$ ranging from $0$ to $\delta$, denote by $\hat{n}_d^+ (m)$ 
(resp. $\hat{n}_d^- (m)$) the total number of these curves which are real, of mass $m$ and
with a non-isolated (resp. isolated) real double point at $y_{c_1 (X) d - 2}$. Define then :
$$\theta_s^d (y,J) = \sum_{m=0}^\delta (-1)^m (\hat{n}_d^+ (m) - \hat{n}_d^- (m)).$$
\begin{theo}
\label{theotheta}
Let $(X, \omega , c_X)$ be a real rational symplectic $4$-manifold, and 
$d \in H_2 (X ; \Z)$.
Let $y \subset X$ be a real configuration of $c_1 (X) d - 2$
distinct points and $s \neq 0$ be the cardinality of $y \cap \R X$. Finally, let 
$J \in \R {\cal J}_\omega$ be an almost complex structure generic enough, so that the
integer $\theta_s^d (y,J)$ is well defined. Then, this
integer $\theta_s^d (y,J)$ neither depends on the choice of $J$ nor on the choice
of $y$ (provided the cardinality of $y \cap \R X$ is $s$).
\end{theo}

For convenience, this integer $\theta_s^d (y,J)$ will be denoted by $\theta_s^d$,
and we put $\theta_s^d = 0$ when $s$ does not have the same parity as $c_1 (X) d$. This 
invariant
makes it possible to give relations between the coefficients of the polynomial
$\chi^d$, namely :
\begin{theo}
\label{theorel}
Let $(X, \omega , c_X)$ be a real rational symplectic $4$-manifold, $d \in H_2 (X ; \Z)$ 
and $r$ be an integer between $0$ and $c_1 (X)d -3$. Then
$\chi^d_{r+2} = \chi^d_r + 2 \theta^d_{r+1}.$ 
\end{theo}

\subsection{Proof of Theorem \ref{theotheta}}
\label{subsect31}

To begin with, we construct as in \S \ref{sectmoduli} the moduli space ${\cal M}_0^d (y)$
of real rational pseudo-holomorphic maps $u : S \to X$, realizing the homology class $d$,
and mapping the marked points $z_1 , \dots , z_{c_1 (X)d - 2}$ of $S$ to the
corresponding points $y_1 , \dots , y_{c_1 (X)d - 2}$ of $X$ and mapping $z_{c_1 (X)d - 1}$
also to $y_{c_1 (X)d - 2}$. This moduli space is obtained by taking the quotient of the
space of such maps by the group ${\cal D}iff^+ (S , z)$ acting by reparametrization. Since
now $u ( z_{c_1 (X)d - 2}) = u ( z_{c_1 (X)d - 1})$, there is a degree two extension of this
group acting by reparametrization, namely the group of diffeomorphisms of $S$, preserving
the orientation, fixing the points $z_1 , \dots , z_{c_1 (X)d - 3}$, and fixing or
exchanging the points $z_{c_1 (X)d - 2}$ and $z_{c_1 (X)d - 1}$. This degree two extension
induces a $\Z / 2\Z$-action on ${\cal M}_0^d (y)$ which has no fixed point, since its
effect is to exchange the two local branches at the double point $y_{c_1 (X)d - 2}$.
Denote by $\widetilde{{\cal M}}_0^d (y)$ the orbit space of this action. It is a Banach 
manifold of class $C^{l, \alpha}$ equipped with an index zero Fredholm projection 
$\tilde{\pi}$ on ${\cal J}_\omega$. Denote by $\tilde{\pi}_\R : \R \widetilde{{\cal M}}_0^d 
(y) \to \R {\cal J}_\omega$ the restriction of $\tilde{\pi}$. The Theorem  of regular values
\ref{theoreg} applies also in this situation, so that the set of regular values of
$\tilde{\pi}$ intersects $\R {\cal J}_\omega$ in a dense set of the second category.

Let then $J_0, J_1 \in \R {\cal J}_\omega$ be regular values of the projection 
$\tilde{\pi} : \widetilde{{\cal M}}_0^d (y) \to {\cal J}_\omega$ such that
no reducible
$J_0$ or $J_1$-holomorphic curve in the class $d$ passes through $y$ with a double point at
$y_{c_1 (X)d - 2}$. Let
$\gamma : [0,1] \to \R {\cal J}_\omega$ be a path transversal to the projection
$\tilde{\pi}_\R : \R \widetilde{{\cal M}}^{d}_0 (y) \to \R {\cal J}_\omega$,
joining $J_0$ to $J_1$. Hence,
$\R \widetilde{{\cal M}}_\gamma = \tilde{\pi}_\R^{-1} (Im (\gamma))$ is a 
submanifold of dimension one
of $\R \widetilde{{\cal M}}^{d}_0 (y)$, equipped with a projection 
$\tilde{\pi}_\gamma : \R \widetilde{{\cal M}}_\gamma
\to [0,1]$ induced by $\tilde{\pi}_\R$.
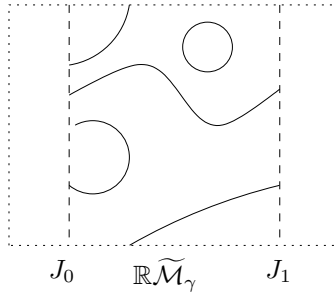
$$\vcenter{\hbox{\input{gro7.pstex_t}}}$$
The path $\gamma$ is chosen so that every element of $\R \widetilde{{\cal M}}_\gamma$ is a 
nodal curve,
with the exception of a finite number of them which may have a unique real ordinary cusp,
a unique real triple point or a unique real tacnode. This path is also chosen so that
when a sequence of elements of $\R \widetilde{{\cal M}}_\gamma$ converges in Gromov 
topology to a
reducible curve of $X$, then this curve has only two irreducible components, both real,
and only nodal points as singularities. Moreover, this path is chosen so that if
$[u , J_S , J] \in \R \widetilde{{\cal M}}_\gamma$ has a unique real ordinary cusp at the 
parameter
$t_u \in \R S$, then the generator $\psi_u$ of $H^0_{D^*} (S , K_S \otimes_\C 
N_{u, -z}^*)_{-1}$ does not vanish at $t_u$. Finally, it is chosen so that
when a sequence of elements of $\R \widetilde{{\cal M}}_\gamma$ converges in Gromov 
topology to an irreducible curve of $X$ not in $\R \widetilde{{\cal M}}_\gamma$, thus 
a cuspidal curve, then this curve has a real
ordinary cusp at $y_{c_1 (X)d - 2}$, and only nodal points as remaining singularities.
Such a choice of $\gamma$ is possible
from Propositions \ref{propcusps}, \ref{proppsiu} and \ref{propprod}.
The integer $\theta_s^d (y, \gamma (\lambda))$ is then well defined for all $\lambda \in 
[0,1]$ but a finite number of values $0 < \lambda_1 < \dots < \lambda_j < 1$ corresponding
either to reducible curves, to cuspidal curves, or to curves having a real triple point
or tacnode. Since the function $\lambda \mapsto \theta_s^d (y, \gamma (\lambda))$ is
obviously constant between these values, we just have to prove that for $i \in \{ 1,
\dots , j \}$, $\theta_s^d (y, \gamma (\lambda_i^-)) = \theta_s^d (y, \gamma (\lambda_i^+))$
where $\lambda_i^-$ (resp. $\lambda_i^+$) is the left limit (resp. right limit) of 
$\lambda$ at $\lambda_i$. The only cases to consider is the apparition of a cuspidal curve,
the cusp being at $y_{c_1 (X)d - 2}$, or the apparition of a curve with a tacnode,
the tacnode being at $y_{c_1 (X)d - 2}$. Indeed, all the other cases follow along the same 
lines as in the proof of Theorem \ref{theoprinc}, the only additional thing to remark is 
that the topology of the node at $y_{c_1 (X)d - 2}$ does not change under these moves. We
will only consider the case of a cuspidal curve, since the other one can be treated exactly 
in the same way. 

So, let $(C_\lambda)_{\lambda \in ]\lambda_i - \epsilon , \lambda_i [}$ be a continuous
family of $\gamma (\lambda)$-holomorphic curves in $\R \widetilde{{\cal M}}_\gamma$ which
converges in Gromov topology to a real cuspidal irreducible $\gamma (\lambda_i)$-holomorphic
 curve, the cusp being at $y_{c_1 (X)d - 2}$. All these curves are nodal as soon as 
$\epsilon$ is small enough. Moreover, such a family is unique. Indeed, if for
$\lambda \in ]\lambda_i - \epsilon , \lambda_i [$, there were two 
$\gamma (\lambda)$-holomorphic curves $C'$ and $C''$ close to the cuspidal curve, then
they would have two intersection points in the neighborhood of each nodal point of
the cuspidal curve, plus four intersection points at $y_{c_1 (X)d - 2}$ and moreover, they
would intersect each other at each point of $y$. This would give
$$(c_1 (X)d - 3) + 4 + 2(\frac{1}{2} (d^2 - c_1 (X)d + 2)-1) = d^2 +1$$
intersection points, which is to much since all multiplicities are positive. In particular,
this family is made of real curves, and since the parity of the number of real curves
does not change, this family does extend to a continuous family 
$(C_\lambda)_{\lambda \in ]\lambda_i - \epsilon , \lambda_i + \epsilon  [}$ of real 
$\gamma (\lambda)$-holomorphic curves. Then, after the transformation, either the topology 
of the real node at $y_{c_1 (X)d - 2}$ is unchanged and then the mass of the curve is
also unchanged, or it has changed, but then the mass of the curve also has changed.
In both cases, the integer $\theta_s^d (y, \gamma (\lambda))$ is left invariant, hence
the result. $\square$

\subsection{Proof of Theorem \ref{theorel}}

Let $y=(y_1, \dots , y_{c_1 (X) d -2})$ be a real configuration of distinct points of $X$,
such that $y_{c_1 (X) d -2} \in \R X$ and $\# ( y \cap \R X) = r+1$. Denote by ${\cal M}_0^d (y)$
the moduli space of rational pseudo-holomorphic curves of $X$ passing through $y$ in the 
homology class $d$. Similarly, denote by $\widetilde{\cal M}_0^d (y)$ the moduli space of such
curves which have a real ordinary node at $y_{c_1 (X) d -2}$. This space has been introduced
in \S \ref{subsect31}. Denote by $P(T_{y_{c_1 (X) d -2}} X)$ the space of tangent lines of $X$ at
$y_{c_1 (X) d -2}$. Then the projection $[u, J_S , J] \in {\cal M}_0^d (y) \mapsto
(J, d|_{z_{c_1 (X) d -2}} u (T_{z_{c_1 (X) d -2}} S) ) \in {\cal J}_\omega \times 
P(T_{y_{c_1 (X) d -2}} X)$ is Fredholm of index zero. Let $(J , \tau) \in \R {\cal J}_\omega 
\times P(T_{y_{c_1 (X) d -2}} \R X)$ be a regular value of this projection. We also assume
that $J$ is a regular value of the projection $\widetilde{\cal M}_0^d (y) \to {\cal J}_\omega$
and that there exists no reducible or cuspidal rational $J$-holomorphic curve passing through
$y$ in the homology class $d$ and having a node or $\tau$ as a tangent at $y_{c_1 (X) d -2}$.
There is then only finitely many element of $\R {\cal M}_0^d (y)$ having $\tau$ as
a tangent at $y_{c_1 (X) d -2}$. These curves are all nodal and irreducible. For every integer
$m$ between $0$ and $\delta$, denote by $\tilde{n}_d (m)$ the number of such curves which are 
real and of mass $m$. Denote then by:
$$\widetilde{\chi}_r^d (y ,J) = \sum_{m=0}^\delta (-1)^m \tilde{n}_d (m).$$
\begin{prop}
\label{proprel}
Under the above assumptions, we have the relations :
$$\chi_{r+2}^d = \widetilde{\chi}_r^d (y ,J) + 2 \sum_{m=0}^\delta (-1)^m \hat{n}_d^+ (m),$$
$$\chi_{r}^d = \widetilde{\chi}_r^d (y ,J) + 2 \sum_{m=0}^\delta (-1)^m \hat{n}_d^- (m).$$
\end{prop}
The integers $\hat{n}_d^+ (m)$ and $\hat{n}_d^- (m)$ have been defined in \S \ref{subsectenonrel}.
The Theorem \ref{theorel} follows easily from this Proposition \ref{proprel} and the definition
of the invariant $\theta$.\\

{\bf Proof of Proposition \ref{proprel}:}

Let us first prove the first relation. For this purpose, let us fix a path 
$\mu : \, ]-\epsilon , \epsilon [ \to \R X$ of class $C^2$ such that $\mu (0) = y_{c_1 (X) d -2}$
and $\mu' (0) \in \tau$. For every $\lambda \in ]-\epsilon , \epsilon [ \setminus \{ 0 \}$,
denote by $y_\lambda$ the set $(y_1, \dots , y_{c_1 (X) d -2}, \mu (\lambda))$. Denote then 
by $\R {\cal M}_0^d (y_\lambda)$ the moduli space of real rational pseudo-holomorphic curves of 
$X$ passing through $y_\lambda$ in the homology class $d$. Then $J$ is a regular value of the 
projection $\pi_\R^{\lambda} : \R {\cal M}_0^d (y_\lambda) \to \R {\cal J}_\omega $ as soon as
$\lambda$ is close enough to zero. Indeed, from Gromov compactness theorem, as soon as
$\lambda$ is close enough to zero, the elements of $\R {\cal M}_0^d (y_\lambda)$ are close,
in Gromov topology, either to elements of $\R \widetilde{\cal M}_0^d (y)$ having a non-isolated
real node at $y_{c_1 (X) d -2}$, or to elements of $\R {\cal M}_0^d (y)$ having $\tau$ as
a tangency at $y_{c_1 (X) d -2}$. As a consequence, these curves are neither cuspidal, nor 
irreducible, and thus $J$ is a regular value of $\pi_\R^{\lambda}$ from Proposition
\ref{proppir}. The set $\{ (\pi_\R^{\lambda})^{-1} (J) , \lambda \in ]-\epsilon , 0 [ \}$ is thus
the union of the images of finitely many continuous functions $C_1 (\lambda), \dots ,
C_j (\lambda)$. Each of these functions converges as $\lambda$ goes to zero either to an
irreducible real $J$-holomorphic curve having a non-isolated real node at $y_{c_1 (X) d -2}$,
or to an irreducible curve having $\tau$ as a tangency at $y_{c_1 (X) d -2}$. We will prove
that each curve of the first kind (resp. second kind) is limit of exactly two (resp. one) such
functions $C_{i_1} (\lambda) , C_{i_2} (\lambda)$. The first relation of Proposition
\ref{proprel} follows, since $\chi_{r+2}^d = \chi_{r+2}^d (y_\lambda ,J)$ for $\lambda$ close
enough to zero, and since the masses of the curves are unchanged while passing to the limit
$\lambda \to 0$.

Let then $C_0$ be an element of $\R {\cal M}_0^d (y)$ having $\tau$ as a tangency at 
$y_{c_1 (X) d -2}$. Since by hypothesis, $(J , \tau)$ is a regular value of the projection
${\cal M}_0^d (y) \to {\cal J}_\omega \times P(T_{y_{c_1 (X) d -2}} X)$, the $J$-holomorphic
curves in a
neighborhood of $C_0$ in $\R {\cal M}_0^d (y) $ are exactly parametrized by their tangencies at
$y_{c_1 (X) d -2}$. Once we move this tangency, we see that these curves provide a foliation
of an angular neighborhood $A$ of $y_{c_1 (X) d -2}$ in $\R X$.
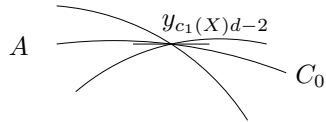
$$\vcenter{\hbox{\input{gro8.pstex_t}}}$$
Since the path $\mu : ]-\epsilon , \epsilon [ \to \R X$ satisfies $\mu (0) = y_{c_1 (X) d -2}$
and $\mu' (0) \in \tau$, restricting $\epsilon$ if necessary, we can assume that its image is
completely included in $A$. For every $\lambda \in ]-\epsilon , \epsilon [ \setminus \{ 0 \}$,
there exists thus one and only one $J$-holomorphic curve in a neighborhood of $C_0$, 
passing through $y_\lambda$, which was the announced result.

Now let $C_0$ be an element of $\R \widetilde{\cal M}_0^d (y)$. This element lifts into two
elements $C_1$ and $C_2$ of the moduli space $\R {\cal M}_0^d (\overline{y})$, where
$\overline{y} = (y_1 , \dots , y_{c_1 (X)d-2} , y_{c_1 (X)d-2})$, see the begining of
\S \ref{subsect31}. Denote by $\R \widehat{\cal M}_0^d (y)$ the moduli space of real rational
pseudo-holomorphic maps having $c_1 (X)d-1$ distincts marked points $z_1 , \dots , z_{c_1 (X)d-1}$
at the source, realizing the homology class $d$ and such that $u(z_i) = y_i$ for 
$1 \leq i \leq c_1 (X)d-2$. The map $\hat{\pi}_\R : [u,J_S ,J] \in \R \widehat{\cal M}_0^d (y)
\mapsto (J , u(z_{c_1 (X)d-1})) \in \R {\cal J}_\omega \times \R X$ is Fredholm of index zero.
The value $(J , y_{c_1 (X)d-2})$ is regular for this projection. Thus the curves $C_1$ and $C_2$
in $(\hat{\pi}_\R)^{-1} (J , y_{c_1 (X)d-2})$ extend in a unique way into two families
$C_1 (\lambda)$ and $C_2 (\lambda)$ of $(\hat{\pi}_\R)^{-1} (J , \mu (\lambda)) = 
({\pi}_\R^\lambda)^{-1} (J)$. These are the two families we were looking for.

The second relation of Proposition \ref{proprel} can be proved in a similar way. We choose this
time a path $\mu : \, ]-\epsilon , \epsilon [ \to X$ of class $C^2$ such that $\mu (0) =
 y_{c_1 (X) d -2}$, $\mu' (0) \in J (\tau)$, and for every $\lambda \in ]-\epsilon , \epsilon [$,
$c_X (\mu (\lambda)) = \mu (- \lambda)$. For every $\lambda \in ]-\epsilon , \epsilon [ 
\setminus \{ 0 \}$,
denote by $y_\lambda$ the set $(y_1, \dots , y_{c_1 (X) d -3}, \mu (\lambda), \mu (-\lambda))$, and
by $\R {\cal M}_0^d (y_\lambda)$ the corresponding moduli space. Now, from Gromov compactness 
theorem, as soon as
$\lambda$ is close enough to zero, the elements of $\R {\cal M}_0^d (y_\lambda)$ are close,
in Gromov topology, either to elements of $\R \widetilde{\cal M}_0^d (y)$ 
having a real isolated node at $y_{c_1 (X) d -2}$ or to elements of $\R {\cal M}_0^d (y)$ 
having $\tau$ as
a tangency at $y_{c_1 (X) d -2}$. Now, each curve of the first kind (resp. second kind) is limit 
of exactly two (resp. one) families $C_{i_1} (\lambda) , C_{i_2} (\lambda)$ of elements of
$\R {\cal M}_0^d (y_\lambda)$. Since the masses of these curves are unchanged while passing to the 
limit $\lambda \to 0$, the second relation follows from the fact that 
$\chi_{r}^d = \chi_{r}^d (y_\lambda ,J)$. $\square$

\subsection{Non-triviality of $\chi^4 (T)$ and $\chi^5 (T)$ for the complex projective plane}

\subsubsection{Generalization of the invariant $\theta$}

The invariant $\theta$ has been defined fixing the position of one of the double points of the 
pseudo-holomorphic curves in the homology class $d$. More generally, one can define such an
invariant fixing the position of $\sigma$ double points of these curves, where 
$0 \leq \sigma \leq \frac{1}{2} [c_1 (X) d - 1]$. More precisely, let
$y = (y_1 , \dots , y_{c_1 (X) d - 1 - \sigma})$ be a real configuration of $c_1 (X) d - 1 - 
\sigma$ distinct points of $X$, and $s$ be the number of those which are real. 
We assume that $y_{c_1 (X) d - 2\sigma}, y_{c_1 (X) d - 2\sigma + 1}, \dots , 
y_{c_1 (X) d - 1 - \sigma}$ are real, so that $s \geq \sigma$. Let 
$J \in \R {\cal J}_\omega$ be generic enough. Then there are 
only finitely many $J$-holomorphic rational curves in $X$ in the 
homology class $d$ passing through $y$ and having a node at each of the points
$y_{c_1 (X) d - 2\sigma}, y_{c_1 (X) d - 2\sigma + 1}, \dots , 
y_{c_1 (X) d - 1 - \sigma}$. These curves are all nodal and 
irreducible. For every integer $m$ ranging from $0$ to $\delta$, denote by $\hat{n}_d^+ (m)$ 
(resp. $\hat{n}_d^- (m)$) the total number of these curves which are real, of mass $m$ and
with an even (resp. odd) number of real isolated double points at 
$y_{c_1 (X) d - 2\sigma}, y_{c_1 (X) d - 2\sigma + 1}, \dots , 
y_{c_1 (X) d - 1 - \sigma}$. Define then :
$$\theta_s^{d,\sigma}  (y,J) = \sum_{m=0}^\delta (-1)^m (\hat{n}_d^+ (m) - \hat{n}_d^- (m)).$$
These definitions extend the ones given in paragraph \ref{subsectenonrel}. In particular,
$\theta_s^{d,0}  (y,J) = \chi_s^d (y,J)$ and $\theta_s^{d,1}  (y,J) = \theta_s^d (y,J)$.
\begin{theo}
\label{theothetagen}
Let $(X, \omega , c_X)$ be a real rational symplectic $4$-manifold,  
$d \in H_2 (X ; \Z)$ and $0 \leq \sigma \leq \frac{1}{2} [c_1 (X) d - 1]$.
Let $y \subset X$ be a real configuration of $c_1 (X) d - 1 - 
\sigma$
distinct points and $s \geq \sigma$ be the cardinality of $y \cap \R X$. Finally, let 
$J \in \R {\cal J}_\omega$ be an almost complex structure generic enough, so that the
integer $\theta_s^{d,\sigma} (y,J)$ is well defined. Then, this
integer $\theta_s^{d,\sigma} (y,J)$ neither depends on the choice of $J$ nor on the choice
of $y$ (provided the cardinality of $y \cap \R X$ is $s$). $\square$
\end{theo}

The proof of this theorem is the same as the one of Theorem \ref{theotheta}. As usual,
this integer $\theta_s^{d,\sigma} (y,J)$ will be denoted by $\theta_s^{d,\sigma}$,
and we put $\theta_s^{d,\sigma} = 0$ when $s$ does not have the suitable parity. 
\begin{theo}
\label{theorelgen}
Let $(X, \omega , c_X)$ be a real rational symplectic $4$-manifold and $d \in H_2 (X ; \Z)$.
Let $\sigma$ be an integer such that $0 \leq 2 \sigma \leq c_1 (X) d - 3$, and $s$ be an integer
between $\sigma$ and $c_1 (X)d -3 - \sigma$. Then
$\theta_{s+2}^{d,\sigma} = \theta_s^{d,\sigma} + 2 \theta_{s+1}^{d,\sigma + 1}. \quad \square$
\end{theo}
The proof of this theorem is the same as the one of Theorem \ref{theorel}.

\subsubsection{Non-triviality of $\chi^4 (T)$ and $\chi^5 (T)$}
\label{subsectnontriv}

In this subparagraph, the real symplectic $4$-manifold $(X, \omega , c_X)$ is the complex
projective plane equipped with its standard symplectic form $\omega_{st}$ and the complex
conjugation $\conj$. We use the canonical identification of $H_2 (\C P^2 ; \Z)$ with $\Z$.
We defined in \S \ref{subsectresults} an invariant $\chi : d \in \Z \mapsto \chi^d (T) \in
\Z [T]$ and have computed it for $d \leq 3$ in Example $1$ of this paragraph.

\begin{prop}
\label{propchi4} 
Let $(X, \omega , c_X)$ be $(\C P^2 , \omega_{st} , \conj)$. Then $\chi^4 (T)$ and $\chi^5 (T)$
are non-zero polynomials of $\Z [T]$.
\end{prop}

\begin{lemma}
\label{lemmacomp}
Let $(X, \omega , c_X)$ be $(\C P^2 , \omega_{st} , \conj)$. Then
$\theta_q^{3} = 1$ for every odd $1 \leq q \leq 7$
$\theta_r^{4,3} = 1$ for every even $4 \leq r \leq 8$ and $\theta_s^{5,6} = 1$ for every even
$6 \leq s \leq 8$.
\end{lemma}

{\bf Proof:}

The proofs are the same in all the cases, so we will prove only the degree $4$ case. Let $y$ be 
a real
configuration of $8$ distincts points in the plane, $r \geq 3$ of which being real. Let $J \in 
{\cal J}_\omega$ be generic enough. There exists then only one $J$-holomorphic rational curve
of degree $4$ in $\C P^2$, passing through $y$, and having its $3$ double points at
$y_6, y_7 , y_8$. Indeed, if there were two of them, they would intersect at each point
$y_1 , \dots , y_5$ with multiplicity at least one, and at each point $y_6, y_7 , y_8$
with multiplicity at least four. This would give an intersection index greater than $16$ which
is impossible. This implies that the corresponding Gromov-Witten invariant is one, since it is
obviously not zero. Now let $J \in \R {\cal J}_\omega$ be generic enough, this unique curve
is real. Denote by $m$ its mass, we have $\theta_r^{4,3} (x , J) = (-1)^m (-1)^m = 1$.
 $\square$ \\

{\bf Proof of Proposition \ref{propchi4} :}

It is a consequence of Theorem \ref{theorelgen} and Lemma \ref{lemmacomp}. $\square$ \\

For instance, the coefficients of the polynomial $\chi^4 (T)$ satisfy the relations \\
$\chi^4_3 = \chi^4_1 + 2 \theta^4_2,$\\
$\chi^4_5 = \chi^4_1 + 4 \theta^4_2 + 4 \theta^{4,2}_3,$\\
$\chi^4_7 = \chi^4_1 + 6 \theta^4_2 + 12 \theta^{4,2}_3 + 8,$\\
$\chi^4_9 = \chi^4_1 + 8 \theta^4_2 + 24 \theta^{4,2}_3 + 32,$ and \\
$\chi^4_{11} = \chi^4_1 + 10 \theta^4_2 + 40 \theta^{4,2}_3 + 80.$ \\
Hence, all these coefficients
cannot vanish simultaneously. Similarly, \\
$\chi^5_{14} = \chi^5_0 + 14 \theta^5_1 + 84 \theta^{5,2}_2 + 280
\theta^{5,3}_3 + 560 \theta^{5,4}_4 + 672 \theta^{5,5}_5 + 448.$

\begin{rem}
It has been observed recently by I. Itenberg, V. Kharlamov and E. Shustin
that the invariant $\chi^d_{3d-1}$ is in fact positive for every $d>0$, thanks to
the research announcement \cite{Mikh} by G. Mikhalkin
\end{rem}

\section*{Appendix}

\appendix

\section{Proof of Theorem \ref{theoreg} in higher genus}
\label{append}

Let $(X , \omega , c_X)$ be a real symplectic $4$-manifold, $g \in \N$, $d \in H_2 (X ; \Z)$
and $x \subset X$ be a real configuration of $c_1 (X)d + (3-n)(g -1) \geq 0$ distincts points. 
Denote by
$\R {\cal M}_g^d (x)_{imm}$ the open subset of $\R {\cal M}_g^d (x)$ made of immersed
pseudo-holomorphic curves and by $\R {\cal M}_g^d (x)_s$ the subspace of curves having a
unique cuspidal point which is real ordinary. The latter is a codimension $n-1$ Banach submanifold,
which can be proved along the same lines as Proposition \ref{propcusps}.

\begin{prop}
\label{prophighgen}
Let $(X , \omega , c_X)$ be a real symplectic manifold of dimension $n$, $g \in \N$, 
$d \in H_2 (X ; \Z)$
and $x \subset X$ be a real configuration of $c_1 (X)d + (3-n)(g -1) \geq 0$ distinct points. 
Then,

1) The space $\{ [u , J_S , J] \in \R {\cal M}_g^d (x)_{imm} \, | \, \dim H^1_D (S , 
N_{u , -z})_{+1} = 1 \}$ is a codimension one Banach submanifold of class $C^{l-1 , \alpha}$ of 
$\R {\cal M}_g^d (x)$ (might be empty).

2) The complementary in $\R {\cal J}_\omega$ of $\pi_\R (\{ [u , J_S , J] \in \R 
{\cal M}_g^d (x)_{imm} \, | \, \dim H^1_D (S , N_{u , -z})_{\pm 1} \geq 1 \})$ is a dense set
of the second category of $\R {\cal J}_\omega$.
\end{prop}

{\bf Proof:}

Let us start with the first part of the proposition. From Lemma $3.2.7$ of \cite{Shev}, the
fibered spaces over $\R {\cal M}_g^d (x)_{imm}$ whose fiber over $[u , J_S , J]$ are the
spaces $L^{k,p} (S , N_{u , -z} )$ and $L^{k-1,p} (S , \Lambda^{0,1} S \otimes_\C N_{u , -z} )$
respectively have the structure of Banach vector bundles of class $C^{l-1 , \alpha}$. Moreover,
the normal Gromov operator $D^N$ induces a $\Z / 2\Z$-equivariant bundle homomorphism 
$L^{k,p} (S , N_{u , -z} ) \to L^{k-1,p} (S , \Lambda^{0,1} S \otimes_\C N_{u , -z} )$. Denote
by $D^N_\R$ the associated morphism $L^{k,p} (S , N_{u , -z} )_{+1} \to L^{k-1,p} 
(S , \Lambda^{0,1} S \otimes_\C N_{u , -z} )_{+1}$ over $\R {\cal M}_g^d (x)_{imm}$. Let
$[u , J_S , J] \in \R {\cal M}_g^d (x)_{imm}$ be such that $\dim H^1_D (S , 
N_{u , -z})_{+1} = 1 $. Since $u$ is immersed and $\ind (D^N_\R) = 0$, it implies that
$\dim H^0_D (S , N_{u , -z})_{+1} = 1 $. From the implicit function theorem, to get the first
part of the proposition, it suffices to prove that the operator:
$$\nabla D^N_\R : T_{[u , J_S , J]} \R {\cal M}_g^d (x) \to Hom ( H^0_D (S , N_{u , -z})_{+1} ,
H^1_D (S , N_{u , -z})_{+1})$$
is surjective. For this purpose, let $\psi_-$ be a generator of $H^0_{D^*} (S , K_S \otimes
 N_{u , -z}^* )_{-1} \cong H^1_D (S , N_{u , -z})_{+1}^*$ and $w_+$ be a generator of
$ H^0_D (S , N_{u , -z})_{+1}$. We are searching for $(v , \stackrel{.}{J}_S , 
\stackrel{.}{J}) \in T_{[u , J_S , J]} \R {\cal M}_g^d (x)$ such that :
\begin{eqnarray}
\label{eqnapp}
\Re e \int_S <\psi_- , \nabla_{(v , \stackrel{.}{J}_S , 
\stackrel{.}{J})}  D^N_\R (w_+)> \neq 0.
\end{eqnarray}
Let us fix $v=0$, $\stackrel{.}{J}_S = 0$, and search for a section
$\stackrel{.}{J} \in L^{k,\alpha} (X , \Lambda^{0,1} X \otimes_\C TX)_{+1}$ which vanishes
along $u(S)$. From formula $(4.2.11)$ of \cite{Shev}, since under these conditions only the
term $[7]$ of this formula is non-zero, the relation (\ref{eqnapp}) becomes :
$$\Re e \int_S <\psi_- , \nabla_{w_+} \stackrel{.}{J} \circ du \circ J_S > \neq 0.$$
Let $U$ be an open subset of $S \setminus z$ small enough such that $u$ restricts to an embedding
from $U$ to $X$, $u (U) \cap u(S \setminus U) = \emptyset$, $c_X (U) \cap U = \emptyset$ and
such that $\psi_-$ and $w_+$ do not vanish on $U$. Let $\alpha$ be a section of the bundle
$\Lambda^{0,1} S \otimes N_{u , -z} $ with support on $U$ such that $\Re e \int_S <\psi_- ,
\alpha > \neq 0$. By integration on the tangent direction to $w_+$, and thus normal to $u(S)$,
we construct a section $\stackrel{.}{J}_1$ of $\Lambda^{0,1} X \otimes_\C TX$ with support in 
a neighborhood of $u(U)$, such that $\nabla_{w_+} \stackrel{.}{J} \circ du \circ J_S =
\alpha$ and $\stackrel{.}{J}_1$ vanishes along $u(S)$. The section $\stackrel{.}{J} =
\stackrel{.}{J}_1 + \overline{c}_X^* \stackrel{.}{J}_1$ is suitable, which proves the first 
part of the proposition.

The same proof leads to the fact that the space $\{ [u , J_S , J] \in \R {\cal M}_g^d (x)_{imm} 
\, | \, \dim H^1_D (S ,  N_{u , -z})_{-1} = 1 \}$ is a Banach submanifold of class 
$C^{l-1 , \alpha}$ of $\R {\cal M}_g^d (x)$ of codimension one (or is empty), which proves the
second part of the proposition in this case. In the general case, let $[u , J_S , J] \in \R 
{\cal M}_g^d (x)_{imm}$ be such that $\dim H^1_D (S ,  N_{u , -z})_{\pm 1} =
\dim H^0_D (S ,  N_{u , -z})_{\pm 1} = h \geq 1$. Denote by $D^N_{\pm 1}$ the operator
$L^{k,p} (S , N_{u , -z} )_{\pm 1} \to L^{k-1,p} 
(S , \Lambda^{0,1} S \otimes_\C N_{u , -z} )_{\pm 1}$. Repeating the same proof as before, we
see that the operator 
$$\nabla D^N_{\pm 1} : T_{[u , J_S , J]} \R {\cal M}_g^d (x) \to Hom ( 
H^0_D (S , N_{u , -z})_{+1} , H^1_D (S , N_{u , -z})_{+1})$$
is non-zero. From the implicit function theorem, there exists then, locally, a submanifold 
$V$ of codimension at least $h^2$ in $\R {\cal M}_g^d (x) \times \R^{h^2 - 1}$ which maps onto 
the subspace $\{ [u , J_S , J] \in \R {\cal M}_g^d (x)_{imm} 
\, | \, \dim H^1_D (S ,  N_{u , -z})_{\pm 1} = h \}$ of $\R {\cal M}_g^d (x)$. The projection
$V \to {\cal J}_\omega$ induced then by $\pi_\R$ is Fredholm of index $-1$, and maps onto
$\pi_\R ( \{ [u , J_S , J] \in \R {\cal M}_g^d (x)_{imm} 
\, | \, \\ \dim H^1_D (S ,  N_{u , -z})_{\pm 1} = h \} )$, hence the second part of the 
proposition. $\square$\\

{\bf Proof of Theorem \ref{theoreg} :}

It is a consequence of Proposition \ref{proptransv} and Proposition \ref{prophighgen}. $\square$

\addcontentsline{toc}{part}{\hspace*{\indentation}Bibliography}

 \nocite{*}  

\bibliography{gromov}

\begin{thebibliography}{10}

\bibitem{Audin}
M.~Audin and J.~Lafontaine, editors.
\newblock {\em Holomorphic curves in symplectic geometry}, volume 117 of {\em
  Progress in Mathematics}.
\newblock Birkh\"auser Verlag, Basel, 1994.

\bibitem{DgKh}
A.~I. Degtyarev and V.~M. Kharlamov.
\newblock Topological properties of real algebraic varieties: {R}okhlin's way.
\newblock {\em Russ. Math. Surv.}, 55(4):735--814, 2000.

\bibitem{Gro}
M.~Gromov.
\newblock Pseudoholomorphic curves in symplectic manifolds.
\newblock {\em Invent. Math.}, 82(2):307--347, 1985.

\bibitem{HLS}
H.~Hofer, V.~Lizan, and J.-C. Sikorav.
\newblock On genericity for holomorphic curves in four-dimensional
  almost-complex manifolds.
\newblock {\em J. Geom. Anal.}, 7(1):149--159, 1997.

\bibitem{IShev}
S.~Ivashkovich and V.~Shevchishin.
\newblock Structure of the moduli space in a neighborhood of a cusp-curve and
  meromorphic hulls.
\newblock {\em Invent. Math.}, 136(3):571--602, 1999.

\bibitem{Kont}
M.~Kontsevich and Y.~Manin.
\newblock Gromov-{W}itten classes, quantum cohomology, and enumerative
  geometry.
\newblock {\em Comm. Math. Phys.}, 164(3):525--562, 1994.

\bibitem{MDSal}
D.~McDuff and D.~Salamon.
\newblock A survey of symplectic {$4$}-manifolds with {$b\sp {+}=1$}.
\newblock {\em Turkish J. Math.}, 20(1):47--60, 1996.

\bibitem{MiWh}
M.~J. Micallef and B.~White.
\newblock The structure of branch points in minimal surfaces and in
  pseudoholomorphic curves.
\newblock {\em Ann. of Math. (2)}, 141(1):35--85, 1995.

\bibitem{Mikh}
G.~Mikhalkin.
\newblock Counting curves via lattice paths.
\newblock {\em Research announcement, math.AG/0209253}, 2003.

\bibitem{RT}
Y.~Ruan and G.~Tian.
\newblock A mathematical theory of quantum cohomology.
\newblock {\em J. Differential Geom.}, 42(2):259--367, 1995.

\bibitem{Shev}
V.~V. Shevchishin.
\newblock Pseudoholomorphic curves and the symplectic isotopy problem.
\newblock {\em preprint math.SG/0010262}, 2000.

\bibitem{SS}
S.~Smale.
\newblock An infinite dimensional version of {S}ard's theorem.
\newblock {\em Amer. J. Math.}, 87:861--866, 1965.

\bibitem{Sot}
F.~Sottile.
\newblock Enumerative real algebraic geometry.
\newblock {\em
  http://www.maths.univ-rennes1.fr/~raag01/surveys/ERAG/index.html}, 2002.
\newblock Electronic survey.

\bibitem{Wels}
J.-Y. Welschinger.
\newblock Invariants of real rational symplectic 4-manifolds and lower bounds
  in real enumerative geometry.
\newblock {\em C. R. Acad. Sci. Paris S\'er. I Math.}, 336(4):341--344, 2003.

\bibitem{Wit}
E.~Witten.
\newblock Two-dimensional gravity and intersection theory on moduli space.
\newblock In {\em Surveys in differential geometry (Cambridge, MA, 1990)},
  pages 243--310. Lehigh Univ., Bethlehem, PA, 1991.

\end{thebibliography}
\bibliographystyle{abbrv}

\noindent Ecole Normale Sup\'erieure de Lyon\\
Unit\'e de Math\'ematiques Pures et Appliqu\'ees\\
$46$, all\'ee d'Italie\\
$69364$, Lyon C\'edex $07$\\
(FRANCE)\\
e-mail : {\tt jwelschi@umpa.ens-lyon.fr}

\end{document}